\documentclass[12pt,a4paper,draft]{article}
\usepackage[english]{babel}
\usepackage{amssymb, amsmath, latexsym, amsfonts, amsthm}
\begin{document}
\begin{center}
ON A SPACE OF INFINITELY DIFFERENTIABLE FUNCTIONS ON AN UNBOUNDED CONVEX SET IN ${\mathbb R}^n$ ADMITTING HOLOMORPHIC EXTENSION IN ${\mathbb C}^n$ AND ITS DUAL 
\footnote
{This work was supported by the grants RFBR 08-01-00779, 08-01-97023 and LSS-3081.2008.1.}
\end{center}
\begin{center}
I.Kh. MUSIN, P.V. FEDOTOVA
\end{center}
\newtheorem{lemma}{Lemma}
\renewcommand{\abstractname}{}
\begin{abstract}
{\sc Abstract}. We consider a space of infinitely smooth functions on an unbounded closed convex set in ${\mathbb R}^n$. 
It is shown that each function of this space can be extended to an entire function in ${\mathbb C}^n$ 
satisfying some prescribed growth condition.  Description of linear continuous  
functionals on this space in terms of their Fourier-Laplace transform is obtained. 
Also a variant of the Paley-Wiener-Schwartz theorem for tempered distributions is given it the paper. 
\end{abstract}

\vspace {1cm}
 
\begin{center}
{\bf \S 1. Introduction} 
\end{center}

{\bf 1.1. On a problem}. 
Let $C$ be an open convex acute cone in ${\mathbb R}^n$ with the apex in the origin [1, p. 73] and  
$b$ be a convex continuous positively homogeneous function of degree 1 on ${\overline C}$ -- the closure of $C$ in ${\mathbb R}^n$. 
The pair $(b, C)$ determines a closed convex unbounded set
$$
U(b, C) = \{\xi \in {\mathbb R}^n: \ 
 -<\xi, y> \ \le b(y),  
\forall y \in C \},
$$  
not containing a straight line. 

Note that the interior of the set $U(b, C)$ is not empty and coincides with the set 
$$
V(b, C)= \{\xi \in {\mathbb R}^n: \ 
-<\xi, y> \ 
 < b(y),  
\forall y \in {\overline C}\}, 
$$
and closure of $V(b, C)$ in ${\mathbb R}^n$ is $U(b, C)$. 

For brevity denote $U(b, C)$ by $U$ and $V(b, C)$ by $V$. 

Let $M=(M_k)_{k=0}^{\infty} $ be a non-decreasing sequence of positive numbers $M_k$ such that $M_0=1$ and 
$$
\displaystyle\lim_{k\rightarrow\infty}\frac {\ln M_k} {k} = +\infty.
$$

Define the space $G_{M}(U)$ of infinitely differentiable functions $f$ on $U$ as follows. For every $m\in\mathbb N$ and $\varepsilon > 0$ let $G_{m,\varepsilon}(U)$ be the space  of $C^{\infty}(U)$-functions $f$ with a finite norm
$$ 
p_{m, \varepsilon}(f)= \sup_{x \in V, \alpha \in {\mathbb Z_+^n}} \frac {\vert (D^{\alpha}f)(x)\vert (1+\Vert x \Vert)^m}{\varepsilon^{\vert \alpha \vert} M_{\vert \alpha \vert}} \ ,  \ m \in {\mathbb N}, \varepsilon > 0.
$$ 
Put $G_M(U)=\displaystyle \bigcap_{m=1}^{\infty}\bigcap_{\varepsilon > 0} G_{m, \varepsilon}(U)$. Thus, $G_M(U)$ is a subclass of the Schwartz class of fast decreasing functions on $U$. 

With usual operations of addition and multiplications by complex numbers $G_M(U)$ is a linear space. 
The family of norms $p_{m, \varepsilon}$ defines a local convex topology in $G_M(U)$.

Note that if $(\varepsilon_m)_{m=1}^{\infty}$ is an arbitrary decreasing to zero sequence of positive numbers $\varepsilon_m$ then the topology in $G_M(U)$ can be determined also by the system of norms
$$ 
p_m(f)= \sup_{x \in V, \alpha \in {\mathbb Z_+^n}} \frac {\vert (D^{\alpha}f)(x)\vert (1+\Vert x \Vert)^m}{\varepsilon_m^{\vert \alpha \vert} M_{\vert \alpha \vert}} \ , \ m \in {\mathbb N}.
$$ 

We are interested to have a description  of the strong dual space to the spaces $G_M (U)$ in terms of the Fourier-Laplace transform 
of continuous linear functionals on $G_M (U)$. 

Detailed consideration of this problem depends on additional conditions on a sequence $M$.  
J.W. de Roever [2] studied this problem under the following assumptions on  $M$: 

1). $M_k^2 \le M_{k-1}M_{k+1},  \ \ \forall k \in {\mathbb N}$; 

2). $ \exists H_1 > 1 \ \exists H_2 > 1$  \ $\forall k, m \in {\mathbb Z_+}$ \ \  
$
M_{k+m} \le H_1 H_2^{k+m} M_k M_m;
$

3). $\exists A >0$ \ $\forall m \in {\mathbb N}$ \ 
$
\displaystyle \sum_{k=m+1}^{\infty}\frac {M_{k-1}}{M_k} \le A m \frac {M_m}{M_{m+1}} \ .
$ 

In this case the sequence $M$ is not quasianalytic. Also from conditions 1) and 3) it follows that there exist numbers $h_1, h_2>0$ such that 
$$
M_k \ge h_1 h_2^k k!, \  
\forall k \in {\mathbb Z_+}.
$$

The same problem was considered in [3](using other methods) under more weak restrictions on the sequence $M$. Namely, the conditions 2) and 3) were replaced with the following conditions:

$2)'$. $ \exists H_1 > 1 \ \exists H_2 > 1$ \  $\forall k \in {\mathbb Z_+}$  \
$
M_{k+1} \le H_1 H_2^k M_k;
$
      
$3)'$.  $ \exists Q_1 >0 \ \exists Q_2>0$ \  $\forall k \in {\mathbb Z_+}$  \
$
M_k \ge Q_1 Q_2^k k!.
$

Here we study this problem assuming that the sequence $M$ with $M_0=1$ satisfies the following conditions:

$i_1)$. $M_k^2 \le M_{k-1}M_{k+1},  \ \ \forall k \in {\mathbb N}$; 

$i_2)$. $ \exists H_1 > 1 \ \exists H_2 > 1$ \ $\forall k, m \in {\mathbb Z_+}$ \ \  $M_{k+m} \le H_1 H_2^{k+m} M_k M_m$;

$i_3)$. $ \forall \varepsilon > 0$ $ \exists a_{\varepsilon} > 0 $ $\forall k \in {\mathbb Z_+}$ \
$M_k \le a_{\varepsilon}\varepsilon^k k!$;

$i_4)$. $ \exists \gamma \in (0, 1) \ \exists b_1 >0 \ \exists  b_2 > 0 \ \forall k \in {\mathbb Z_+}$ \ $M_k \ge b_1 b_2^k k!^{\gamma}$;

$i_5)$. there exists a logarithmically convex non-decreasing sequence 
$K=(K_m)_{m=0}^{\infty}$ with $K_0=1$ such that for some $t_1>1, t_2 > 1$ 
$$
t_1^{-1} t_2^{-m} K_m \le \frac {m!}{M_m} \le t_1 t_2^m K_m, \ m \in {\mathbb Z_+}.
$$

The sequence $(m!^{\alpha})_{m=0}^{\infty}$ (where $\alpha \in (0, 1)$) is a simple example of such a sequence.

{\bf 1.2. Definitions and notations}. For $u=(u_1, \ldots , u_m) \in {\mathbb R}^m ({\mathbb C}^m)$, $v=(v_1, \ldots , v_m) \in {\mathbb R}^m ({\mathbb C}^m)$ 
let $<u, v> = u_1v_1 + \cdots + u_m v_m$, $\Vert u \Vert$ be the Euclidean norm in ${\mathbb R}^m ({\mathbb C}^m)$. 

For $z\in\mathbb C^m$, $R>0$ let $B_R (z)$ be a ball in $\mathbb C^m$ of a radius $R$ with the center at the point $z$. Let $\nu_m (R) = \nu_m(1)R^{2m}$ be a volume of $B_R(z)$. 

Denote the tubular domain ${\mathbb R^n} + i C$ by $T_C$. 

If $\Omega \subset {\mathbb R}^m$ ($\Omega \subset {\mathbb C}^m $) then the distance from 
$x \in {\mathbb R}^m$ ($z \in  {\mathbb C}^m $) to the set $\Omega$ we denote by $d_{\Omega}(x) (d_{\Omega}(z))$. 
The distance from  $x \in \Omega$ ($z \in  \Omega$) to the boundary of $\Omega$ is denoted by $\Delta_{\Omega}(x) (\Delta_{\Omega}(z))$.

For a locally convex space $X$ let $X'$ be the space of linear continuous functionals on $X$ and let $X^*$ be the strong dual space. 

For open set $\Omega$ in ${\mathbb C}^m$ \ $H(\Omega)$ is a set of holomorphic functions in $\Omega$ , $psh(\Omega)$ is a set of plurisubharmonic  functions in $\Omega$.

With a sequence $L=(L_k)_{k=0}^{\infty} $  of positive numbers 
$L_k $ with $L_0=1$ such that 
$$
\displaystyle \lim  \limits_{k \to \infty} \frac {\ln L_k}{k} = + \infty
$$ 
let us associate the function
$\omega_L$:
$
\omega_L(r) = \displaystyle \sup_{k \in {\mathbb Z_+}} \ln\frac {r^k} {L_k}$ for $r>0$ and $\omega_L(0) = 0$. 

If decreasing to zero sequence $(\varepsilon_m)_{m=1}^{\infty}$ of numbers $\varepsilon_m$ is chosen then for brevity denote $\omega_M(\frac{r}{\varepsilon_m})$ by $\omega_m(r), \ r \ge 0$. 

{\bf 1.3. Main results}. Note that from the condition $i_3)$ on the sequence $M$ it follows that each function belonging to $G_M(U)$ admits holomorphic continuation in ${\mathbb C}^n$. Using the conditions $i_4)$ and $i_5)$ it can be shown (see theorem 1) that the space  $G_M(U)$ is topologically isomorphic to the space $E(U)$ of entire functions $f$ in ${\mathbb C}^n$ such that for every $\varepsilon > 0, m \in {\mathbb N}$ there exists a constant $C_{m, \varepsilon} > 0$ such that
$$
\vert f(z) \vert \le C_{m, \varepsilon} \frac {e^{w_K(\varepsilon d_U(x))  + w_K(\varepsilon \Vert y \Vert)}}{(1 + \Vert z \Vert)^m} \ , \ z \in {\mathbb C}^n,
$$
endowed with a topology defined by the family of norms 
$$
q_{m, \varepsilon}(f) = \sup_{z \in {\mathbb C}^n} 
\frac {\vert f(z) \vert (1 + \Vert z \Vert)^m}{e^{w_K(\varepsilon d_U(x))  + w_K(\varepsilon \Vert y \Vert)}} \ , \ \varepsilon > 0, m \in {\mathbb N}.
$$
Here as usual $x = Re  z, y = Im  z$. 

For description of dual space to $G_M(U)$ we use a method of paper [4]. But to use it first we have to study the Fourier-Laplace transform of tempered distributions with support in $U(b, C)$. That's why we consider the space $S(U)$ of $C^{\infty}(U)$-functions $f$ 
such that for each $p \in {\mathbb Z_+}$
$$
{\Vert f \Vert}_{p, U} = \sup_{x \in V, \vert \alpha \vert \leq p} \vert (D^\alpha f)(x) \vert (1 + {\Vert x \Vert})^p  < \infty. 
$$
Let $S_p(U)$ be a completion of $S(U)$ by the norm ${\Vert \cdot \Vert}_{p, U}$.  
We endow $S(U)$ with a topology of projective limit of the spaces $S_p(U)$. 
It is known that $S^*(U)$ is topologically isomorphic to the space of tempered distributions with support in $U(b, C)$ [5].

For each $m \in {\mathbb N}$ let
$$
V_{b, m}(T_C) = \{f \in H(T_C): N_m(f) = \sup_{z\in T_C} 
\frac
{\vert f(z)\vert e^{-b(y)}} 
{(1 + \Vert z \Vert)^m(1 + \frac {1} {\Delta_C(y)})^m} < \infty\}, 
$$
$$
H_m(T_C) = \{ F \in H(T_C): {\Vert F \Vert}_m = \sup_{z \in T_C} 
\frac 
{\vert F(z) \vert}
{e^{\omega_m (\Vert z \Vert)} (1 + \frac {1} {\Delta_C (y)})^m} < \infty \} \ , 
$$
where $z = x +i y, x \in {\mathbb R}^n, y \in C$. 

Put $H_M(T_C) = \bigcup_{m=1}^\infty H_m(T_C)$, $V_b(T_C) = \bigcup_{m=1}^\infty V_{b, m}(T_C)$.
With usual operations of addition and multiplication by complex numbers $H_M(T_C)$ and $V_b(T_C)$ are linear spaces. 
Supply $H_M(T_C)$ ($V_b(T_C)$) with the topology of inductive limit of spaces $H_m(T_C)$ ($V_{b, m}(T_C)$). 

{\it The Fourier-Laplace transform of a functional} $\varPhi \in S^*(U)$ ($\varPhi \in G_M^*(U)$) is defined by formula 
$$\hat \varPhi(z) = (\varPhi, e^{i<\xi, z>}), \ z \in T_C.$$ 

{\it The Laplace transform of a functional} $\varPhi \in E^*(U)$  is defined by formula 
$$\tilde \varPhi(z) = (\varPhi, e^{i<\lambda, z>}), \ z \in T_C.$$ 

The following results on description of dual spaces are obtained in the paper.

{\bf Theorem 2.} {\it The Fourier-Laplace transform ${\cal F}: S^*(U) \to V_b(T_C)$  establishes a topological isomorphism between the  spaces $S^*(U)$ and $V_b(T_C)$.}

If $b(y) = a \Vert y \Vert ( a \ge 0)$ then we have well-known result of V.S. Vladimirov [1]. The main part in the proof of theorem 2 is proving that   ${\cal F}$ is acting "into". 

{\bf Theorem 3.} {\it The Fourier-Laplace transform establishes a topological isomorphism between the  spaces $G_M^*(U)$ and $H_M(T_C)$}.  

{\bf Theorem 4.} {\it The Laplace transform establishes a topological isomorphism between the spaces $E^*(U)$ and $H_M(T_C)$}. 

The proof of theorem 1 is given in the third section. The proof of theorem 2 is  in the forth section. Theorem 3 is proved in the fifth section and  theorem 4 is proved in the sixth section. In the second section some properties of functions  associated with the sequences $M$ and $K$ are given. A simple example of the set $U(b, C)$ is given in [3].

{\bf Remark 1}. The definition of the spaces $G_M(U)$ and $H_{b, M}(T_C)$ does not depend on a choice of the sequence $(\varepsilon_m)_{m=1}^{\infty}$. So we put ${\varepsilon}_m = H_2^{-m}, \ 
m \in {\mathbb N}$. Thus, $\omega_m(r)= \omega_M(H_2^m r), \ r \ge 0$.

\begin{center}
{\bf \S 2. Auxiliary results} 
\end{center}
\begin{lemma} 
For $x>0$ and $s>1$ 
$$
x \le \ln \frac {s}{s-1} + \sup_{m \in {\mathbb Z}_+} \ln \frac {(sx)^m}{m!} \ .
$$
\end{lemma} 

{\bf Proof}. Let $x>0, s>1$. We have
$$
e^x = 
\sum_{m = 0}^{\infty} 
\frac {(sx)^m}{m!} 
\frac {1} {s^m} 
\le \sup_{m \in {\mathbb Z}_+} 
\frac {(sx)^m}{m!}
\sum_{m =0}^{\infty}  \frac {1}{s^m} = \frac {s}{s-1} \sup_{m \in {\mathbb Z}_+} 
\frac {(sx)^m}{m!} \ .
$$
From this the assertion follows.

\begin{lemma} 
Let $ s > 1 $. Then for $r \ge 0$
$$
w_M(sr) + r \le w_M(sH_2r) + \ln \frac {a_1 H_1 s}{s-1} \ .
$$
\end{lemma}

{\bf Proof}. The lemma is true for $r=0$. Now let $r > 0 $. Then for $ s > 1$ 
$$
w_M(sr) + r \le 
\sup_{k \in {\mathbb Z}_+} 
\ln \frac {(sr)^k}{M_k} + 
\ln \frac {s}{s-1} + \sup_{m \in {\mathbb Z}_+} 
\ln \frac {(sr)^m}{m!} \le 
$$
$$ 
\le
\sup_{k \in {\mathbb Z}_+} \ln \frac {(sr)^k}{M_k} + \ln \frac {s}{s-1} + 
\sup_{m \in {\mathbb Z}_+} \ln \frac {a_1(sr)^m}{M_m}
=
$$
$$
= \sup_{k, m \in {\mathbb Z}_+} \ln \frac {(sr)^{k+m}}{M_k M_m} + \ln \frac {a_1 s}{s-1} \le
$$
$$
\le \sup_{k, m \in {\mathbb Z}_+} \ln \frac {H_1(sH_2r)^{k+m}}{M_{k + m}} + \ln \frac {a_1 s}{s-1} 
\le w_M(sH_2r) + 
\ln \frac {a_1 H_1 s}{s-1} \ .
$$

Lemma 2 is proved.

Applying lemma 2 with $s=H_2^m$ $( m \in {\mathbb N}$) we have for each $ r \ge 0$
$$
w_M(H_2^m r) + r \le  w_M(H_2^{m+1}r) + \ln \frac {a_1 H_1 H_2^m}{H_2^m-1} \le  w_M(H_2^{m+1}r) + \ln \frac {a_1 H_1 H_2}{H_2-1} \ . 
$$
Setting $Q=\displaystyle \ln \frac {a_1 H_1 H_2}{H_2-1}$ and  taking into consideration definitions of $w_m$ we have for each $m \in {\mathbb N}$
$$
w_m(r) + r \le  w_{m+1}(r) + Q, \ r \ge 0.  \eqno (1)
$$

Define a function ${\cal N}$ on $[0, \infty)$ as follows:
$$
{\cal N}(r)=\min \left \{k \in {\mathbb Z_+}: w_M(r) = 
\displaystyle \ln 
\frac {r^k}{M_k} \right \}, \ r>0; {\cal N}(0) = 0.
$$ 

It is easy to check that for
$r \in (\frac {M_k}{M_{k-1}}, \frac {M_{k+1}}{M_k}]$ $(k \in {\mathbb N})$ \ 
$w_M(r)= \displaystyle \ln \frac {r^k}{M_k}$,  ${\cal N}(r) = k$; $w(r) = 0, {\cal N}(r)=0$ if $r \in (0, M_1]$. 
It is clear that $\omega_M$ is continuous on $[0, \infty)$.  

From the condition $i_4)$ on the sequence $M$ it follows that there exists a constant $A_{\gamma}>0$ such that 
$$
w_M(r) \le A_{\gamma} r^{\frac 1 \gamma}, \ r \ge 0. \eqno (2)
$$

Using the condition $i_3)$ on $M$ and lemma 1 we have for each $\varepsilon \in (0, 1)$
$$
w_M(r) \ge \sup_{k \in {\mathbb Z}_+} 
\ln \frac {r^k}{a_{\varepsilon^2}\varepsilon^{2k} k!} \ge \frac {r}{\varepsilon} - \ln \frac {a_{\varepsilon^2}}{1 - \varepsilon} \ , \ r>0.
$$

Thus,
$
 \displaystyle \lim_{r \to + \infty} \frac {w_M(r)}{r} = + \infty.
$

Now we shall estimate the growth of  ${\cal N}(r)$. 
Using the representation 
$$
w_M(r) = \int_0^r \frac {{\cal N}(t)}{t} \ dt, \  r \ge 0, \eqno (3)
$$
we obtain
$$
w_M(er)\ge \int_r^{er} \frac {{\cal N}(t)}{t} \ dt \ge {\cal N}(r).
$$
Now using (2) we have
$$
{\cal N}(r) \le A_{\gamma} (e r)^{\frac 1 \gamma}, \ r \ge 0. \eqno (4)
$$

\begin{lemma} 
For each $r_1, r_2 \ge 0$
$$
\vert w_M(r_2) -  w_M(r_1) \vert \le A_{\gamma} e^{\frac 1 \gamma}(r_1 +r_2)^{\frac 1 \gamma - 1} \vert r_2 - r_1 \vert.
$$
\end{lemma}

{\bf Proof}. Let $r_2 \ge r_1 \ge 0$. 
Using (3) and (4) we have
$$
w_M(r_2) -  w_M(r_1) = \int_{r_1}^{r_2} \frac {{\cal N}(t)}{t} \ dt \le \gamma A_{\gamma} e^{\frac 1 \gamma} ({r_2}^{\frac 1 \gamma} - {r_1}^{\frac 1 \gamma}) \le
$$
$$
\le A_{\gamma} e^{\frac 1 \gamma}(r_1 + r_2)^{\frac 1 \gamma - 1} (r_2 - r_1).
$$

Lemma 3 is proved.

\begin{lemma}
For each $N \in {\mathbb N}$ 
$$
w_K(r) + N \ln r \le w_K(e r) + \ln  K_N, \ r > 0.
$$
\end{lemma}

{\bf Proof}. Since $M$ is a logarithmically convex sequence then for  $p, q \in {\mathbb Z_+}$ 
\  $M_{p+q} \ge M_p M_q$. So
$$
K_{p+q} = \frac {(p+q )!}{M_{p+q}} \le e^{p+q} \frac {p! q!}{M_pM_q}
= e^{p+q}K_p K_q \ , p, q \in {\mathbb Z_+}.
$$ 
Let $r > 0$. Then for each $N \in {\mathbb N}$ 
$$
w_K(r)+ N \ln r \le  \sup_{m \in {\mathbb Z_+}} \ln \frac {r^m}{K_m} + \ln r^N = 
\sup_{m \in {\mathbb Z_+}} \ln \frac {r^{m+N}}{K_m} \le 
$$
$$
\le 
\sup_{m \in {\mathbb Z_+}} 
\ln \frac {(er)^{m+N}K_N}{K_{m+N}} 
\le w_K(e r) + \ln K_N.
$$
Lemma 4 is proved.

By lemma 4 for $N \in {\mathbb N}$ and $A>1$ 
$$
w_K(r) + N \ln (1 + Ar) \le w_K(e r) + N \ln 2A + \ln K_N, \ r \ge 0. \eqno (5)
$$

\begin{lemma} For each $r \ge 0$ \ 
$
2w_K(r) \le w_K(e r). 
$
\end{lemma}

{\bf Proof}. For $r=0$ it is obvious. For $r>0$ we have
$$
2w_K(r)= 2 \sup_{m \in {\mathbb Z_+}} \ln \frac {r^m}{K_m} = 
\sup_{m \in {\mathbb Z_+}} \ln \frac {r^{2m}}{K_m^2} \leq
$$
$$
\leq 
\sup_{m \in {\mathbb Z_+}} \ln \frac {(er)^{2m}}{K_{2m}} 
\le w_K(e r).
$$

\begin{center}
{\bf \S 3. Isomorphism of spaces $G_M(U)$ and $E(U)$ }. 
\end{center}

{\bf Theorem 1.} {\it The spaces $G_M(U)$ and $E(U)$ are topologically isomorphic}.

{\bf Proof}. 
Let $f \in G_M(U)$ be an  arbitrary. Then 
$\forall \alpha \in {\mathbb Z_+^n} \ \forall  \varepsilon > 0 \ \forall m \in {\mathbb N}$
$$
\vert (D^{\alpha}f)(x)\vert \le p_{m, \varepsilon}(f) 
\frac 
{\varepsilon^{\vert \alpha \vert} M_{\vert \alpha \vert}} {(1+\Vert x \Vert)^m} \ , \ x \in V.  \eqno (6)
$$
From (6) it follows that for points $x, x_0 \in V$ 
$$
f(x)=\displaystyle \sum_{\vert \alpha \vert \ge 0} \frac {(D^{\alpha}f)(x_0)}{\alpha!} (x-x_0)^{\alpha}
$$
and the series standing to the right of this equality converges uniformly on the compact subsets of $V$ to $f$.

Construct an isomorphism $T: G_M(U)\rightarrow E(U)$ as follows. Let $x_0 \in V$. In view  of (6) and the condition $i_3)$ on $M$  
$$
F_{x_0}(z)= \displaystyle \sum_{\vert \alpha \vert \ge 0} \frac {(D^{\alpha}f)(x_0)}{\alpha!} (z-x_0)^{\alpha} \ , \ z \in {\mathbb C}^n,
$$
is an entire function.
Note that for $x \in V$ \ $F_{x_0}(x)=f(x)$. So for $x_1, x_2 \in V$ \ 
$F_{x_1}(z)= F_{x_2}(z), \ z \in {\mathbb C}^n$. 
Thus, we defined the function $F \in H({\mathbb C}^n)$ such that for each $\xi \in V$ we have  $F= F_{\xi}$ in ${\mathbb C}^n$ 
 and $F(x)=f(x), \ x \in V$. Put $T(f)=F$. Obviously, the mapping $T$ is one-to-one and linear.  

Let $f \in G_M(U)$. Now we shall estimate the growth of $F=T(f)$. 
Let $z = x +~iy, \ x \in V, y \in {\mathbb R}^n$. 
Since 
$$
F(z) = \displaystyle \sum_{\vert \alpha \vert \ge 0} \frac 
{(D^{\alpha}f)(x)}{\alpha!} (iy)^{\alpha} \ , 
$$
then for each $m \in {\mathbb N}$ and $\varepsilon > 0$
$$ 
\vert F(z) \vert \le p_{m, \varepsilon}(f)
\sum_{\vert \alpha \vert \ge 0} \frac 
{\varepsilon^{\vert \alpha \vert} M_{\vert \alpha \vert} \Vert y \Vert^{\alpha}} {(1+\Vert x \Vert)^m \alpha!} = 
$$
$$
=
\frac {p_{m, \varepsilon}(f)} {(1+\Vert x \Vert)^m} 
\sum_{N=0}^{\infty}
\varepsilon^N M_N \Vert y \Vert^N 
\sum_{\vert \alpha \vert = N} \frac {1}{\alpha!} = 
$$
$$
= \frac {p_{m, \varepsilon}(f)} {(1+\Vert x \Vert)^m}
\sum_{N=0}^{\infty}
\varepsilon^N M_N \Vert y \Vert^N \frac {n^N}{N!} \le \frac {2 t_1 p_{m, \varepsilon}(f)} {(1+\Vert x \Vert)^m}
\sup_{N \in {\mathbb Z_+}}
\frac 
{(2\varepsilon n t_2 \Vert y \Vert)^N} {K_N} =  
$$
$$
=\frac {2 t_1p_{m, \varepsilon}(f)} {(1+\Vert x \Vert)^m} e^{w_K(2\varepsilon n t_2 \Vert y \Vert)}. \eqno (7)
$$
Now we estimate  $\vert F(z) \vert$ at the points $z = x +iy$ such that $ x \notin V, y \in {\mathbb R}^n $. 
Let $\xi$ be an arbitrary point of  $V$. Then from the representation 
$$
F(z)= \displaystyle \sum_{\vert \alpha \vert \ge 0} \frac {(D^{\alpha}f)(\xi)}{\alpha!} (z-\xi)^{\alpha}
$$
we have for each $m \in {\mathbb N}$ and $\varepsilon > 0$ 
$$
\vert F(z) \vert \le p_{m, \varepsilon}(f) 
\sum_{\vert \alpha \vert \ge 0} 
\frac 
{\varepsilon^{\vert \alpha \vert} M_{\vert \alpha \vert}}
{(1 + \Vert \xi \Vert)^m \alpha!} \Vert z-\xi \Vert^{\vert\alpha \vert} = 
$$
$$
= \frac {p_{m, \varepsilon}(f)} {(1 + \Vert \xi \Vert)^m}
\sum_{N=0}^{\infty} 
\frac
{\varepsilon^N M_N \Vert z-\xi \Vert^N n^N}{N!}
 \le 
$$
$$
\le 
\frac {t_1 p_{m, \varepsilon}(f)} {(1 + \Vert \xi \Vert)^m}
\sum_{N=0}^{\infty} 
\frac
{(\varepsilon n t_2 \Vert z-\xi \Vert)^N}{K_N} \le
$$
$$
\le 
\frac {2 t_1 p_{m, \varepsilon}(f)} {(1 + \Vert \xi \Vert)^m} e^{w_K(2\varepsilon n t_2 \Vert z - \xi \Vert)}. 
$$
Thus, in this case for each $m \in {\mathbb N}$ and $\varepsilon > 0$
$$
\vert F(z) \vert \le 
2 t_1 p_{m, \varepsilon}(f) \inf_{\xi \in V}
\frac {e^{w_K(2\varepsilon n t_2 \Vert z - \xi \Vert)}} {(1 + \Vert \xi \Vert)^m}. \eqno (8)
$$
For $m \in {\mathbb N}$ and $s > 0$ let 
$$
g_{m, s}(z) = \inf \limits_{\xi \in V}(w_K(s \Vert z - \xi \Vert) - m \ln (1 + \Vert \xi \Vert)), 
$$
where $z=x+iy$ and $x\notin V$, $y\in\mathbb R^n$. 
 
Since
$$
w_K(s \Vert z - \xi \Vert) - m \ln (1 + \Vert \xi \Vert) 
\leq
$$
$$
\le
w_K(2 s \Vert x - \xi \Vert) - m \ln (1 + \Vert \xi \Vert)+ 
w_K(2s\Vert y \Vert) \leq
$$
$$
\leq 
w_K(2s \Vert x - \xi \Vert) + m \ln (1 + \Vert x - \xi \Vert)+ 
w_K(2s \Vert y \Vert)- m \ln (1 + \Vert x \Vert), 
$$
then
$$
g_{m, s}(z) \le w_K(2sd_U(x)) + m \ln (1 + d_U(x))+ w_K(2s \Vert y \Vert) - m \ln (1 + \Vert x \Vert).
$$
Using inequality $(5)$ and putting $d_{m,s}=m\ln \left(1+\frac {1}{2s}\right) + \ln R_m$ we have for $z=x+i y$ with $x \notin V$, $y\in\mathbb R^n$
$$
g_{m,s}(z)\leq \omega_K (2esd_U(x)) + \omega_K(2s\Vert2\Vert) - m\ln(1+\Vert x\Vert) + d_{m,s}. \eqno (9)
$$
Going back to (8) and using (9) with $s=2\varepsilon  n t_2$ we obtain
$$
\vert F(z) \vert \leq 
A_{m, \varepsilon} 
p_{m, \varepsilon}(f) 
e^{w_K(4\varepsilon e n t_2 d_U(x))  + w_K(4\varepsilon  n t_2 \Vert y \Vert) - m \ln (1 + \Vert x \Vert)}, \eqno (10)
$$
where 
$ x \notin V, y \in {\mathbb R}^n $,
$A_{m, \varepsilon}= 2 t_1 e^{d_{m, s}}$. Using (7) we conclude that (10) holds everywhere in ${\mathbb C}^n$.

Now we shall continue the estimate (10). For each $z \in {\mathbb C}^n$
$$
\vert F(z) \vert \le 
A_{m, \varepsilon} 
p_{m, \varepsilon}(f) 
e^{w_K(4\varepsilon e n t_2 d_U(x))  + w_K(4\varepsilon  n t_2 \Vert y \Vert) - m \ln (1 + \Vert z \Vert) + m \ln (1 + \Vert y \Vert) }. 
$$
Again using inequality (5) and putting $B_{m,\varepsilon}=K_m A_{m,\varepsilon}(1+\frac {1} {4\varepsilon nt_2})^m$  we have for each $\varepsilon >0, m\in\mathbb N$
$$
\vert F(z) \vert \leq
B_{m, \varepsilon} 
p_{m, \varepsilon}(f) 
e^{w_K(4\varepsilon e n t_2 d_U(x))  + w_K(4\varepsilon  n t_2 \Vert y \Vert) - m \ln (1 + \Vert z \Vert)}, \ z \in \mathbb C^n.
$$  
Thus, for each $\varepsilon >0, m\in\mathbb N$
$$
q_{m,4e\varepsilon nt_2}(T(f))\leq B_{m,\varepsilon}p_{m,\varepsilon}(f).
$$
This means that  $T$ is a continuous mapping from $G_M(U)$ to $E(U)$. 

Now we prove that the inverse mapping $T^{-1}$ is continuous. 
Let $F \in E(U)$. Show that $f=F_{|U}$ belongs to $G_M(U)$. 
Let $m \in {\mathbb N}$, $\varepsilon \in (0, 1), R>0$ are arbitrary. 
Let $x \in V$. For every $\alpha \in {\mathbb Z_+^n}$
$$
(1 + \Vert x \Vert)^m (D^{\alpha}f)(x) = 
\frac {\alpha! }{(2\pi i)^n} 
\displaystyle {\int \cdots \int}_{L_R(x)}
\frac 
{(1 + \Vert x \Vert)^m F(\zeta)}
{(\zeta_1 - x_1)^{\alpha_1 +1} \cdots (\zeta_n - x_n)^{\alpha_n +1}} \ d \zeta ,
$$
where  
$L_R(x)= \{\zeta = (\zeta_1, \ldots , \zeta_n) \in {\mathbb C}^n: \vert \zeta_j - x_j \vert = R, j=1, \ldots , n \}$, 
$d \zeta = d \zeta_1 \cdots d  \zeta_n $.
From this
$$
 (1 + \Vert x \Vert)^m \vert (D^{\alpha} f)(x) \vert 
\le \frac 
{\alpha!}{R^{\vert \alpha \vert}}  
\max_{\zeta \in L_R} (1 + \Vert  \zeta - x \Vert)^m (1 + \Vert  \zeta \Vert)^m \vert F(\zeta) \vert
\le 
$$
$$
\le 
\frac 
{\alpha!}{R^{\vert \alpha \vert}}(1 + \sqrt n R)^m
q_{m, \varepsilon}(F)
\max_{\zeta = \xi + i \eta \in L_R}
e^{w_K(\varepsilon d_U(\xi))  + w_K(\varepsilon \Vert \eta \Vert)} \le 
$$
$$
\le
\frac 
{\alpha!}
{R^{\vert \alpha \vert}}
q_{m, \varepsilon}(F)
e^{2w_K(\varepsilon \sqrt n R) + m \ln (1 + \sqrt n R)}.  
$$
Using lemma 5 we have
$$
(1 + \Vert x \Vert)^m \vert (D^{\alpha} f)(x) \vert 
\le  \frac 
{\vert \alpha \vert!}
{R^{\vert \alpha \vert}}
q_{m, \varepsilon}(F)
e^{w_K(e\varepsilon \sqrt n R) + m \ln (1 + \sqrt n R)} .
$$
Now using inequality (5) we have for each $R>0$
$$
(1 + \Vert x \Vert)^m \vert (D^{\alpha} f)(x) \vert \le
\frac 
{\vert \alpha \vert!}
{R^{\vert \alpha \vert}}K_m
q_{m, \varepsilon}(F) 
e^{w_K(e^2\varepsilon \sqrt n R) + m \ln \frac {2}{e\varepsilon}} .
$$
Consequently,
$$
(1 + \Vert x \Vert)^m \vert (D^{\alpha} f)(x) \vert \le \frac {\vert \alpha \vert!}{K_{\vert \alpha \vert}} K_m
(e^2\varepsilon \sqrt n)^{\vert \alpha \vert}
\left(1+\frac {1}{e\varepsilon}\right)^m q_{m, \varepsilon}(F)
\le
$$
$$
\le 
t_1 (e^2r_2 \varepsilon \sqrt n)^{\vert \alpha \vert} M_{\vert \alpha \vert} K_m
\left(1+\frac {1}{e\varepsilon}\right)^m q_{m, \varepsilon}(F).
$$
This means that
$$
p_{m, e^2t_2 \varepsilon \sqrt n}(f) \le t_1  K_m  
\left(1+\frac {1}{e\varepsilon}\right)^m q_{m, \varepsilon}(F).
$$
Thus, the mapping $T^{-1}$ is continuous. 

Theorem 1 is proved.

\begin{center}
{\bf \S 4. Description of $S^*(U)$ in  terms of the Fourier-Laplace transform} 
\end{center}

Let $C^*=\{\xi\in{\mathbb R}^m: \  <\xi,x> \ \geq 0,  \forall x \in C \} $ be a dual cone to $C$, 
$pr \ C$ be an intersection of $C$ with a unique sphere. 
For $r \ge 0$ let $B_r=\{\xi \in {\mathbb R}^n: \Vert \xi \Vert \le r \}$,  $\widetilde B_r$ be an exterior of $B_r$.

\begin{lemma} 
For  $y \in C, m \in {\mathbb N}$ let 
$$
g(\xi) = -<\xi, y> + m \ln(1 + \Vert\xi\Vert), \ \xi \in {\mathbb R}^n.
$$ 
Then there exists a number $d>0$ not depending on $y$ such that
$$
\sup \limits_{\xi \in U}g(\xi)
\le 
b(y) + d m + 3m \ln \left(1 + \frac {1}{\Delta_C(y)}\right) + 2m \ln (1 + \Vert y \Vert).
$$
\end{lemma}

{\bf Proof}. Since $b$ is  continuous and positively homogeneous on ${\overline C}$ then there exists a number $r > e $ 
such that $\vert b(y) \vert \le r \Vert y \Vert$ for $y \in {\overline C}$. 
So $U \subset C^* +~B_r$. 
There is a number $R_0>0$ such that for all $R > R_0$ the set $U_R=U \cap B_R$ is not empty. 
Let $\widetilde U_R= U \setminus U_R$, $I_R = (C^* +~B_r)\cap B_R$, $\widetilde I_R = (C^* +~B_r) \setminus I_R$.

Let $\xi_0 \in V$ be  arbitrary. We  will show that for each $y \in C$ there exists a number $R_1>R_0$ such that 
$\displaystyle \sup_{\xi \in \widetilde I_{R_1}} g(\xi) < g(\xi_0)$.  
Note that if $R > R_2=\max(R_0, 2r + \frac {m}{\Delta_C(y)})$
then
$$
\displaystyle \sup_{\xi \in \widetilde T_R} g(\xi) \le 
\displaystyle \sup_{\xi_1 \in C^* \cap {\widetilde B_{R-r}}} g(\xi_1) + 
\displaystyle \sup_{\xi_2 \in B_r} g(\xi_2) =
$$
$$
=\displaystyle \sup_{\sigma \in pr \ C^*} \displaystyle \sup_{t > R-r} (-t <\sigma, y> + m \ln (1 +t)) + \displaystyle \sup_{\xi_2 \in B_r} g(\xi_2)=
$$
$$
\le \displaystyle \sup_{\sigma \in pr \ C^*}(-(R-r) <\sigma, y>) + m \ln (1 + R - r) + r \Vert y \Vert + m \ln (1 +r) =
$$
$$
= -(R-r) \Delta_C (y) + m \ln (1 + R - r) + r \Vert y \Vert + m \ln (1 +r). \eqno (11)
$$
Show that for each $y \in C$ there exists a number $R_3 > 0$ such that for $R \ge R_3$
the following inequality holds
$$
-(R-r) \Delta_C (y) + (2m+ (r + \Vert \xi_0 \Vert) \Vert y \Vert) \ln (1 + R - r) < 0.  \eqno (12)
$$
Then we will have
$$
-(R-r) \Delta_C (y) + m \ln (1 + R - r) + r \Vert y \Vert + m \ln (1 +r) < g(\xi_0). \eqno (13)
$$
Note that the set of solutions of an inequality $x - \lambda \ln (1+x) > 0$ with a parameter $\lambda > 1$ contains the interval $[\lambda^2, \infty)$.
Let $\lambda=\frac {2m+ (r + \Vert \xi_0 \Vert) \Vert y \Vert}{\Delta_C (y)} \ . $ 
Since $r > 1$ then $\lambda > 1$. 
So the inequality (12) holds for all
$R \ge R_3 = r + \left(\frac {2m+ (r + \Vert \xi_0 \Vert) \Vert y \Vert}{\Delta_C (y)}\right)^2$. 
Let 
$R_1 = 3r + R_0 + \frac {m}{\Delta_C (y)} +  \left(\frac {2m + (r + \Vert \xi_0 \Vert) \Vert y \Vert}{\Delta_C (y)}\right)^2 $.
Then from (11) and (13) it follows that $\displaystyle \sup_{\xi \in \widetilde U_{R_1}} g(\xi) < g(\xi_0)$. 
This means that
$\sup \limits_{\xi \in U} g(\xi) = \sup \limits_{\xi \in U_{R_1}} g(\xi)$. 
So the point $\tilde \xi$ at which the upper bound of function $g$ on $U$ is attained belongs to $U_{R_1}$. 
Making elementary estimates we find a constant $d>0$ depending on $r, R_0, m$ and $\xi_0$ such that
$$
\ln (1 + \Vert \tilde \xi \Vert) \le d + 3 \ln \left(1 + \frac {1}{\Delta_C (y)}\right) + 2 \ln (1 + \Vert y \Vert).
$$
Hence 
$$
\sup \limits_{\xi \in U} g(\xi) = g(\tilde \xi) = - <\tilde \xi, y> + m \ln (1 + \Vert \tilde \xi \Vert) \le 
$$
$$
\le 
b(y) + d m + 3m \ln \left(1 + \frac {1}{\Delta_C(y)}\right) + 2m \ln (1 + \Vert y \Vert). 
$$

Lemma 6 is proved.

{\bf Proof of theorem 2}. 
Let $\varPhi \in S^*(U)$ be an arbitrary functional. Then  $\hat \varPhi$ is a holomorphic function in 
$T_C$ ([1], [2], [5]-[8]). 
So there exist numbers $m\in {\mathbb N}$ and $c>0$ such that
$$
\vert (\varPhi, f) \vert \leq c {\Vert f \Vert}_{m, U} , \ f \in S(U). 
$$ 
Putting $f(\xi) = e^{i<\xi ,z>}$ (here $z = x + iy, x \in {\mathbb R^n}, y \in C$) we have
$$
\vert \hat \varPhi(z) \vert \leq c \sup_{\xi \in V, \vert \alpha \vert \le m } {\vert (iz)^\alpha e^{i<\xi, z>} \vert(1 + \Vert \xi \Vert)^m} \le 
$$ 
$$
\le  c (1 + \Vert z \Vert)^m 
e^{\sup \limits_{\xi \in V} (-<\xi, y> + m \ln(1 + \Vert \xi \Vert)}. 
$$ 
Using lemma 6 we obtain
$$
\vert \hat \varPhi(z) \vert \le  c (1 + \Vert z \Vert)^{3m} 
e^{ b(y) + d m} \left(1 + \frac {1}{\Delta_C(y)}\right)^{3m}, \ z \in T_C. 
$$
So the mapping ${\cal F}$ acts from $S^*(U)$ to $V_b(T_C)$.  

It is known [5, p. 20] that $S^*(U)$ is an inductive limit of increasing sequence of Banach spaces $S_p^*(U)$. 
So if 
$\varPhi \in S_m^*(U)$ then 
$$
\vert (\varPhi, f) \vert \leq {\Vert \varPhi \Vert}_{-m, U} {\Vert f \Vert}_{m, U}, \ f \in S(U). 
$$  
Here ${\Vert \varPhi \Vert}_{-m, U}$ is a norm of $\varPhi$ in $S_m^*(U)$, $m \in {\mathbb N}$. 
Now using lemma 6 we have
$$
N_{3m}(\hat \varPhi) \le  e^{d m} {\Vert \varPhi \Vert}_{-m, U}.  
$$
This means that ${\cal F}$ is continuous.

Proving that ${\cal F}$ is bijective and ${\cal F}^{-1}$ is continuous is the same as in [1]. 
Thus, ${\cal F}$ is an isomorphism.

{\bf Remark 2}. 
Let $\varepsilon > 0$, $\eta \in pr C$. At the end of paper [8] it was shown by J.W. de Roever that if $\varPhi \in S^*(U)$ then there exist numbers $C_{\varepsilon} > 0$ and $m \in {\mathbb N}$ 
(not depending on $\varepsilon > 0$) such that 
$$
\vert \hat \varPhi(z) \vert \le C_{\varepsilon}(1 + \Vert z \Vert)^m e^{b(y)}, \ z = x + i y, y \in \varepsilon \eta + C.
$$ 

\begin{center}
{\bf \S 5. Space $G_M(U)$ and its dual}. 
\end{center}

Using the Arzela-Askoli theorem it is easy to show that $G_M(U)$ is the  $(M^*)$-space (definition of  $(M^*)$-spaces see in [9], [10]).

Let
$$
C_m(U)=\{f\in C(U): \widetilde{p}_m (f) = 
\sup_{x \in U} \vert f(x) \vert (1+\Vert x \Vert)^m < \infty \}, \  m \in {\mathbb N}.
$$
By standard scheme ([4, Propositions 2.10, 2.11, Corollary 2.12]) one can prove the following lemma. 

\begin{lemma}
Let $T\in G_M'(U)$ and numbers $c>0$, $m\in\mathbb N$ are such that
$$
\vert (T, f ) \vert \leq c p_m(f), \ f \in G_M(U).
$$
 
Then there exist functionals $T_{\alpha}\in C_m'(U)$ ($\alpha\in\mathbb Z_+^n$) such that
$$
\vert (T_{\alpha},f) \vert\leq\frac{C\widetilde{p}_m(f)}  {\varepsilon_m^{\vert\alpha\vert}M_{\vert\alpha\vert}} \ , \ f \in C_m(U),
$$ 
and 
$
(T, f)=\displaystyle\sum_{\vert\alpha\vert\geq 0}(T_{\alpha}, D^{\alpha}f ), \ 
f \in G_M(U). 
$
\end{lemma}

\begin{lemma} 
Let $S \in G'_M(U)$. Then $\hat S \in H_M(T_C)$. 
\end{lemma}

{\bf Proof}. First note that for each $z=x + i y \in T_C$ ($x \in {\mathbb R}^n, y \in C$) 
the function $f_z(\xi) = e^{i<\xi ,z>}$ belongs to  $G_M(U)$. Indeed, 
for each $m \in {\mathbb N}$ 
$$
p_m(f_z) = 
\sup_{\xi \in V, \alpha \in {\mathbb Z_+^n}} \frac {\vert (iz)^\alpha e^{i<\xi, z>} \vert(1 + \Vert \xi \Vert)^m} {\varepsilon_m^{\vert \alpha \vert} M_{\vert \alpha \vert}} \le
$$
$$
\le 
\sup_{\alpha \in {\mathbb Z_+^n}} 
\frac {{\Vert z \Vert}^{\vert \alpha \vert}}  {\varepsilon_m^{\vert \alpha \vert} M_{\vert \alpha \vert} }
\sup \limits_{\xi \in V}
e^{ -<\xi, y> + m \ln(1 + \Vert \xi \Vert)} =
e^{\omega_m(\Vert z \Vert) +\sup \limits_{\xi \in V}( -<\xi, y> + m \ln(1 + \Vert \xi \Vert))}. 
$$
Using lemma 6 and inequality (1) we can find a constant $A>0$ depending only on $m$ such that
$$
p_m(f_z) \le A e^{\omega_{m + [r] + 1} (\Vert z \Vert)} \left(1 + \frac {1}{\Delta_C(y)}\right)^{3m} \ . \eqno (14)
$$
Here $r$ is a number which was defined in the proof of lemma 6. 

Now let $S \in G'_M(U)$. It is clear that the function $\hat S(z) = (S, e^{i<\xi, z>})$
is correctly defined on $T_C$. 
Using lemmas 7 and 6, condition $i_4)$ it is easy to show that $\hat S \in H(T_C)$. 

Since there exist numbers $m \in {\mathbb N}$ and $c>0$ such that
$$
\vert (S, f) \vert \leq c p_m(f), \ \ f\in G_M(U). 
$$
then using (14) we obtain
$$
\vert \hat S(z) \vert \leq  cA e^{\omega_{m + [r] + 1} (\Vert z \Vert)} \left(1 + \frac {1}{\Delta_C(y)}\right)^{3m} \ .  
$$
Thus, $\hat S \in H_M(T_C)$.

Obviously for each $m \in {\mathbb N}$ 
embeddings $j_m : H_{b, m}(T_C) \rightarrow H_{b, m+1}(T_C)$ are completely continuous. 
So $H_M(T_C)$ is $(LN^*)$-space 
(definition of  $(LN^*)$-spaces see in [9], [10]).

In lemmas 9 and 10 the following notations will be used. 
A point 
$z = (z_1, \ldots , z_k) \in {\mathbb C}^k$ ($k= 2, 3, \ldots $) will be 
written in the form  $z = (z', z_k)$, where $z'= (z_1, \ldots , z_{k-1}) \in {\mathbb C}^{k-1}$. 

For $k=1,\ldots, n-1$ and $z = (z_1, \ldots , z_k) \in {\mathbb C}^k$, $\zeta \in {\mathbb C}^n$
let
$\varphi_k(z,\zeta)=
\varphi(z_1+\zeta_1,\ldots,z_k+\zeta_k,\zeta_{k+1},\ldots,\zeta_n)$, 
for $z = (z_1, \ldots , z_n) \in {\mathbb C}^n$ and $\zeta \in {\mathbb C}^n$ let 
$\varphi_n(z,\zeta) = 
\varphi(z_1+\zeta_1, \ldots, z_n + \zeta_n).
$

If  $f = \sum_{k=1}^m f_k d{\overline z}_k$ is a form of type (0, 1) in ${\mathbb C}^m$ then put 
$\Vert f(z) \Vert^2 = \sum_{k=1}^m \vert f_k(z) \vert^2$. 

For $u\in C^1 (\Omega)$ ($\Omega$ is an open set in $\mathbb C^m$) $\overline\partial u =\sum_{j=1}^m\frac{\partial u} {\partial \overline z_j} d \overline z_j$. For a definition of the  $\overline \partial$ operator on forms see [11].

By $\lambda_m$ denote the Lebesgue measure in ${\mathbb C}^m$.

\begin{lemma} 
Let $\cal O$ be a domain of holomorphy in $\mathbb C^n$. Let $\varphi\in psh (\mathbb C^n)$ and for some
$c_{\varphi}>0$ and $\nu>0$  
$$
\vert\varphi(z)-\varphi(t)\vert\leq c_{\varphi}
$$
if
$\Vert z-t\Vert\leq\frac{1} {(1+\Vert t
\Vert)^{\nu}}$. Let $h\in psh(\cal O)$ and for some $c_h>0$
$$
\vert h(z) - h(t) \vert\leq c_h
$$
if
$\Vert z - t \Vert\leq \min \left(1,\frac{\Delta_{\cal O}(t)} {4}\right)$. 
 
Let $f\in H(\mathbb C^{k-1}\times\cal O)$ (here $k=2,\ldots,n$) and 
 for some $c_f>0, m \geq 0$   
$$
\vert
f(z', \zeta)\vert\leq
c_f(1+ \Vert (z',\zeta)\Vert)^m \left(1 + \frac{1}
{\Delta_{\cal O}(\zeta)}\right)^m
e^{\varphi_{k-1}(z', \zeta) + h(\zeta)}.
$$ 
Then there exists a holomorphic function $F$ in $\mathbb C^k \times \cal O$ such that 
$F(z', 0, \zeta)=f(z', \zeta)$ for 
$z' \in {\mathbb C}^{k-1}, \zeta \in \cal O$ and for some $C>0$, $m_k\geq 0$ 
$$
\vert
F(z, \zeta)\vert\leq
C(1+\Vert(z,\zeta)\Vert)^{m_k} \left(1 + \frac{1}
{\Delta_{\cal O}(\zeta)}\right)^{m_k}
e^{\varphi_k(z,\zeta) + h(\zeta)}, \  z \in \mathbb C^k, \zeta \in \cal O.
$$ 
\end{lemma}

{\bf Proof}.
Let the function $\mu\in C^{\infty}[0,\infty)$ be such that for $t \geq 0$ \ $0\leq \mu(t)\leq 1$ and $\vert \mu'(t) \vert \leq 4$, 
$\mu (t)=1$ for $t\in [0,\frac {1}{3}], \mu(t)=0$ for $t\geq 1$.
For $z \in {\mathbb C}^k$, $\zeta \in {\mathbb C}^n$ let
$$
H_k (z,\zeta)=\mu ((1+\sqrt 2
\Vert(z',\zeta)\Vert)^{\nu}\vert z_k\vert).
$$ 
Note that 
$H_k(z,\zeta)=0$ out of the set 
$$
\Omega_k =\{(z,\zeta)\in\mathbb
C^{k + n}: (1+ \sqrt 2 \Vert (z', \zeta) \Vert)^{\nu} 
\vert z_k \vert < 1 \}.
$$

For $j=1, \ldots , n$
$$\frac {\partial H_k}
{\partial\overline{\zeta}_j}(z,\zeta) = \frac {\sqrt 2}{2} \nu \vert z_k \vert 
\mu'((1+\sqrt 2 \Vert(z', \zeta)\Vert)^{\nu}\vert z_k \vert)(1+\sqrt 2\Vert(z',\zeta)\Vert)^{\nu -1} \frac {\zeta_j} {\Vert (z',\zeta)\Vert}. 
$$
For every $j=1,\ldots, k-1$ 
$$
\frac {\partial H_k}
{\partial\overline{z}_j}(z,\zeta) 
= \frac {\sqrt 2}{2} \nu \vert z_k \vert
\mu'((1+\sqrt 2 \Vert(z',\zeta)\Vert)^{\nu}\vert z_k\vert)(1+\sqrt 2\Vert(z',\zeta)\Vert)^{\nu -1} \frac {z_j} {\Vert (z',\zeta)\Vert}. 
$$
Further
$$\frac {\partial H_k}{\partial\overline{z}_k}(z,\zeta)= 
\frac {1}{2}\mu^{'} ((1+\sqrt 2 \Vert(z^{'},\zeta)\Vert)^{\nu}\vert z_k\vert)(1+\sqrt 2\Vert(z^{'},\zeta)\Vert)^{\nu}\frac {z_k} {\vert z_k \vert}. 
$$
Let 
$$
W_k= \{(z,\zeta)\in\mathbb C^{n+k}: \frac
{1}{3}<(1+\sqrt 2\Vert(z^{'},\zeta)\Vert)^{\nu}\vert z_k \vert<1\}.
$$
Obviously, $\overline{\partial}H_k (z,\zeta)=0$ out of $W_k$. 
If
$(z,\zeta)\in W_k$, then
$$
\Vert\overline{\partial} H_k (z,\zeta)\Vert^2 = 
\sum_{j=1}^k \left\vert \frac {\partial H_k} {\partial {\overline z}_j}(z,\zeta)\right\vert^2 + \sum_{j=1}^n \left\vert \frac {\partial H_k} {\partial \overline{\zeta}_j}(z,\zeta)\right\vert^2 =
$$
$$
=\frac {(\mu'(1+\sqrt 2 \Vert(z',\zeta)\Vert^{\nu} \vert z_k \vert) 
(1+\sqrt 2\Vert(z', \zeta)\Vert)^{\nu})^2}{4} \left(\frac {2\nu^2 \vert z_k \vert^2} 
{(1+\sqrt 2\Vert(z',\zeta)\Vert)^2} + 1\right)\leq
$$
$$\leq 4 (1+\sqrt 2\Vert(z',\zeta)\Vert)^{2\nu} \left(\frac {2\nu^2} {(1+\sqrt 2\Vert(z', \zeta)\Vert)^{2+2\nu}}+1\right)\leq 
$$
$$
\le 
4(2\nu^2 +1)(1+\sqrt 2 \Vert (z',\zeta)\Vert)^{2\nu}.
$$

Thus, everywhere in $\mathbb C^k \times \cal O$ 
$$
\Vert
\overline{\partial} H_k (z,\zeta)\Vert ^2 \leq 4 (2\nu^2 + 1)(1 +
\sqrt 2\Vert (z',\zeta)\Vert)^{2\nu}.
$$
Choose a function $v_k\in C^{\infty}(\mathbb C^k\times\cal O)$ with suitable bound so that the function
$$
F(z,\zeta)=f(z', \zeta)H_k(z,\zeta)-z_k v_k(z,\zeta), \ (z,\zeta)\in\mathbb C^k\times\cal O,
$$ 
(for which $F(z' ,0,\zeta) = f(z', \zeta)$ if $(z',\zeta)\in\mathbb C^{k-1}\times\cal O$) is holomorphic in $\mathbb C^k\times\cal O$.
So $v_k\in C^{\infty}(\mathbb C^k\times\cal O)$ must satisfy the equation 
$$
\overline{\partial} v_k (z,\zeta) = 
\frac {f(z', \zeta )\overline{\partial} H_k (z,\zeta)} {z_k}, \ 
(z,\zeta) \in \mathbb
C^k\times\cal O \eqno (15)
$$
Denote by $g_k(z,\zeta)$ the form of a type $(0,1)$ standing to the right of (15). 
Note that $g_k (z,\zeta)$ is a zero form out of $W_k$. 
Directly can be checked that
$\overline{\partial}
g_k (z,\zeta)=0$ in $\mathbb C^k\times\cal O$. 
If $(z,\zeta)\in
W_k$, then 
$$
\Vert g_k(z,\zeta)\Vert ^2 = 
\frac {\vert
f(z',\zeta)\vert ^2 \Vert \overline{\partial} H_k
(z,\zeta)\Vert ^2} {\vert z_k\vert ^2} 
\leq 
$$
$$
\le 36(2\nu^2+1)
c_f^2 4^{\nu}
(1 + \Vert(z',\zeta)\Vert)^{2m_{k-1}+4 \nu} 
e^{2(\varphi_{k-1}(z',\zeta)+
h(\zeta) + m_{k-1} \ln (1+\frac {1}{\Delta_{\cal O}(\zeta)}))}.
$$
Putting 
$A_1=36(2\nu^2+1)c_f^2 4^{\nu}$
we have in ${\mathbb C^k} \times \cal O$
$$
\Vert g_k(z,\zeta)\Vert ^2\leq A_1 (1 +
\Vert(z',\zeta)\Vert)^{2m_{k-1}+4\nu}e^{2(\varphi_{k-1}(z',\zeta)+
h(\zeta) + m_{k-1} \ln (1+\frac {1}{\Delta_{\cal O}(\zeta)}))}.
$$
Note that for $(z,\zeta)\in\Omega_k$ 
$$
\Vert
(z_1+\zeta_1,\ldots, z_k+\zeta_k,\zeta_{k+1},\ldots,\zeta_n)-
(z_1+\zeta_1,\ldots, z_{k-1}+\zeta_{k-1},\zeta_k,\ldots,\zeta_n)\Vert \leq
$$
$$
\le \frac {1} {(1+\sqrt 2\Vert(z', \zeta)\Vert)^{\nu}} 
\le 
\frac
{1}
{(1+\Vert( z_1+\zeta_1, \ldots, z_{k-1}+\zeta_{k-1}, \zeta_k,\ldots, \zeta_n)\Vert)^{\nu}}.
$$ 
Hence for 
$(z,\zeta)\in\Omega_k$ \ 
$
\vert \varphi_k (z,\zeta) - \varphi_{k-1}(z',\zeta)\vert 
\leq c_{\varphi}.
$
Now we can obtain an integral estimate on $\Vert g_k (z, \zeta) \Vert$. We have
$$
\int_{\mathbb C^k\times\cal O} 
\frac 
{\Vert g_k (z,\zeta)\Vert^2 e^{-2(\varphi_k (z, \zeta) + h(\zeta) + m_{k-1} \ln (1 + \frac {1} {\Delta_{\cal O}(\zeta)}))}}
{(1+\Vert(z,\zeta)\Vert)^{2(n+m_{k-1}+k+\nu)-1}} 
d\lambda_{n+k} (z,\zeta)=
$$
$$
= \int_{W_k} 
\frac {\Vert g_k (z,\zeta)\Vert^2 
e^{-2(\varphi_k (z, \zeta) + h(\zeta) + m_{k-1} \ln (1 + \frac {1} {\Delta_{\cal O}(\zeta)}))}}   
{(1+\Vert(z,\zeta)\Vert)^{2(n+m_{k-1}+k+\nu)-1}} d\lambda_{n+k} (z,\zeta)\leq 
$$
$$
\leq A_1 \int_{W_k} 
\frac {(1+\Vert (z', \zeta)\Vert)^{2\nu} e^{2(\varphi_{k-1}(z',\zeta)-\varphi_k (z,\zeta))}} {(1+\Vert(z,\zeta)\Vert)^{2(n+k)-1}} d\lambda_{n+k} (z,\zeta)\leq
$$
$$
\leq A_1 e^{2c_{\varphi}} \int_{W_k} \
\frac {(1+\Vert (z', \zeta)\Vert)^{2\nu}} {(1+\Vert(z,\zeta)\Vert)^{2(n+k)-1}} d\lambda_{n+k} (z,\zeta)\leq
$$
$$
\leq A_1 e^{2c_{\varphi}} \int \limits_{\mathbb C^{n+k-1}} \int \limits_{\vert z_k\vert < \frac {1}{(1+\sqrt 2\Vert(z', \zeta)\Vert)^{\nu}}}
\frac {(1+\Vert (z', \zeta)\Vert)^{2\nu} } {(1+\Vert(z,\zeta)\Vert)^{2(n+k)-1}}  d\lambda_1 (z_k) d\lambda_{n+k-1}(z',\zeta)=
$$
$$
= \pi A_1 e^{2c_{\varphi}}\int_{\mathbb C^{n+k-1}} \frac {d\lambda_{n+k-1}(z',\zeta)} {(1+\Vert (z', \zeta)\Vert)^{2(n+k-1)+1}}<\infty.$$
Theorem $2.2.1'$ in [12] provides a solution 
$v_k\in C^{\infty}(\mathbb C^k \times \cal O)$ of the equation $(15)$ such that
$$
2 \int_{\mathbb C^k \times \cal O} 
\frac 
{\vert v_k (z,\zeta)\vert^2 
e^{-2(\varphi_k (z, \zeta) + h(\zeta) + m_{k-1} \ln (1 + \frac {1}{\Delta_{\cal O}(\zeta)}))}}  
{(1+\Vert(z,\zeta)\Vert)^{2(n+m_{k-1}+k+\nu)-1} 
(1+\Vert (z,\zeta)\Vert^2)^2} 
\ d \lambda_{n+k} (z,\zeta) 
\leq 
$$
$$
\leq \int_{\mathbb C^k \times\cal O} 
\frac {\vert v_k (z,\zeta)\vert^2 
e^{-2(\varphi_k (z, \zeta) + h(\zeta) + m_{k-1} \ln (1 + \frac {1}{\Delta_{\cal O}(\zeta)}))}}
{(1 + \Vert(z, \zeta)\Vert)^{2(n+m_{k-1}+k+\nu)-1}} 
\  d \lambda_{n+k} (z,\zeta).
$$
Further we have  
$$
\int_{\mathbb C^k\times\cal O} \frac {\vert F (z,\zeta)\vert^2 
e^{-2(\varphi_k (z, \zeta) + h(\zeta) + m_{k-1} \ln (1 + \frac {1}{\Delta_{\cal O}(\zeta)}))}}
 {(1+\Vert(z,\zeta)\Vert)^{2(n+m_{k-1}+k+\nu)-1}(1+\Vert(z,\zeta)\Vert^2)^3} 
 \ d\lambda_{n+k} (z,\zeta)
\leq  
$$
$$
\leq 2 \int_{\mathbb C^k\times\cal O} \frac {\vert f (z', \zeta) \vert^2 \vert H_k (z,\zeta)\vert^2 
e^{-2(\varphi_k (z, \zeta) + h(\zeta) + m_{k-1} \ln (1 + \frac {1}{\Delta_{\cal O}(\zeta)}))}}  
{(1+\Vert(z,\zeta)\Vert)^{2(n+m_{k-1}+k+\nu)-1} 
(1+\Vert(z,\zeta)\Vert^2)^3} \ d\lambda_{n+k} (z,\zeta) + 
$$
$$
+ 2\int_{\mathbb C^k\times\cal O} 
\frac {\vert z_k\vert^2 \vert v_k (z,\zeta)\vert^2 
e^{-2(\varphi_k (z, \zeta) + h(\zeta) + m_{k-1} \ln (1 + \frac {1}{\Delta_{\cal O}(\zeta)}))}} 
{(1+\Vert(z,\zeta)\Vert)^{2(n+m_{k-1}+k+\nu)-1}(1+\Vert(z,\zeta)\Vert^2)^3} \ d\lambda_{n+k} (z,\zeta) \leq
$$
$$
\leq 
2 c_f^2 \int_{\Omega_k} \frac {e^{2(\varphi_{k-1}(z', \zeta)-\varphi_k (z, \zeta))} } {(1+\Vert(z,\zeta)\Vert)^{2(n+k+\nu)-1} (1+\Vert(z,\zeta)\Vert^2)^3} \ d \lambda_{n+k} (z,\zeta) + 
$$
$$
+ 2\int_{\mathbb C^k\times\cal O} \frac {\vert v_k (z,\zeta)\vert^2 
e^{-2(\varphi_k (z,\zeta)+h(\zeta)) + m_{k-1} \ln (1 + \frac {1}{\Delta_{\cal O}(\zeta)}))}
d\lambda_{n+k} (z,\zeta)} {(1+\Vert(z,\zeta)\Vert)^{2(n+m_{k-1}+k+\nu)-1}(1+\Vert(z,\zeta)\Vert^2)^2}\leq
$$
$$
\leq 
2 c_f^2 e^{2c_{\varphi}}\int_{\mathbb C^{n+k}}\frac {d\lambda_{n+k} (z,\zeta)} {(1+\Vert(z,\zeta)\Vert^2)^{n+k+\nu+2}}+
$$
$$
+\int_{\mathbb C^k\times\cal O} \frac {\Vert g_k (z,\zeta)\Vert^2 
e^{-2(\varphi_k (z,\zeta)+h(\zeta)) +  m_{k-1} \ln (1 + \frac {1}{\Delta_{\cal O}(\zeta)}))}
d\lambda_{n+k} (z,\zeta)} {(1+\Vert(z,\zeta)\Vert)^{2(n+m_{k-1}+k+\nu)-1}} = A_2<\infty.
$$

Now from integral estimate we obtain uniform estimate on $F$ by standard arguments (see, for example, [13, p. 205]). 
Let
$R=\min(\frac {\Delta_{\cal O}(\zeta)} {4},\frac {1}
{2(1+\sqrt2\Vert(z,\zeta)\Vert)^{\nu}})$. 
We have in $\mathbb C^k\times\cal O$ 
$$
\vert F(z,\zeta)\vert^2\leq\frac {1}
{\nu_{n+k}(R)}\int \limits_{B_R(z,\zeta)}\vert F(t,w)\vert^2 \
d\lambda_{n+k}(t,w) \le 
$$
$$
\leq 
A_2
\sup_{(t,w) \in B_R(z,\zeta)}
(e^{2(\varphi_k (t,w) + h(w)  + m_{k-1} \ln(1+\frac {1} {\Delta_{\cal O}(w)}))} 
(1+\Vert (t,w)\Vert)^{2(n+m_{k-1}+k+\nu)+5}) \eqno (16) 
$$ 
Since $\Vert w-\zeta\Vert \leq\min
\left(1,\frac {\Delta_{\cal O}(\zeta)}{4}\right)
$
then $\vert h(w) - h(\zeta)\vert\leq c_h$. 
Note also that for $(t,w)\in
B_R(z,\zeta)$
$$
\Vert
(t_1+w_1,\ldots,t_k+w_k,w_{k+1},\ldots,w_n) -
(z_1+\zeta_1,\ldots,z_k+\zeta_k,\zeta_{k+1},\ldots,\zeta_n)\Vert =
$$
$$
=
\Vert (t_1,\ldots,t_k,0,\ldots,0) + w -
(z_1,\ldots,z_k,0,\ldots,0) - \zeta \Vert \leq \Vert t-z\Vert +
\Vert w-\zeta \Vert \leq 
$$
$$
\le 2R \le 
\frac {1} {(1+\sqrt 2 \Vert
(z,\zeta)\Vert)^{\nu}} \leq \frac {1} {(1 +
\Vert(z_1+\zeta_1,\ldots,z_k+\zeta_k,\zeta_{k+1},\ldots,\zeta_n)\Vert)^{\nu}}.
$$
Hence for $(t,w)\in B_R(z,\zeta)$ \ 
$
\vert \varphi_k (t,w)-\varphi_k (z,\zeta)\vert 
\leq c_{\varphi}. 
$ 
From this and $(16)$ we have for $(z,\zeta)\in \mathbb C^k \times \cal O$
$$
\vert F(z,\zeta)\vert^2 \leq  \frac {A_2}{\nu_{n+k}(1)} (2 + \Vert(z,\zeta)\Vert)^{2(n+m_{k-1}+k+\nu)+5} 
\left(1+\frac {2} {\Delta_{\cal O}(\zeta)}\right)^{2m_{k-1}} \cdot
$$
$$
\cdot\left(\frac {4} {\Delta_{\cal O}(\zeta)} + 2(1 + \sqrt2\Vert(z,\zeta)\Vert)^{\nu}\right)^{2(n+k)}
 e^{2(c_{\varphi}+c_h)}e^{\varphi_k (z, \zeta)+h(\zeta)}.
$$
Thus, there exists a constant $C>0$ such that  in  $\mathbb C^k \times
\cal O$ 
$$
\vert F(z,\zeta)\vert \leq
C(1+\Vert(z,\zeta)\Vert)^{m_k}\left(1+\frac {1} {d_{\cal O}(\zeta)}\right)^{m_k} e^{\varphi_k(z,\zeta)+h(\zeta)},
$$
where
$m_k=(n+k)(\nu +1) + m_{k-1} + 3$.

Lemma is proved.

\begin{lemma} 
Let $\cal O$ be a domain of holomorphy in $\mathbb C^n$. 
Let $\varphi\in psh (\mathbb C^n)$ be such that for some
$c_{\varphi}>0$ and $\nu>0$ 
$
\vert\varphi(z)-\varphi(t)\vert\leq c_{\varphi}
$
if $\Vert z-t\Vert\leq\frac{1} {(1+\Vert t
\Vert)^{\nu}}$,  $h\in psh(\cal O)$ be such that for some $c_h>0$
$
\vert h(z) - h(t) \vert\leq c_h
$
if $z, t \in\cal O$ and $\Vert z - t \Vert\leq \min
\left(1,\frac{\Delta_{\cal O}(t)} {4}\right)$. 

Let $f \in H(\cal O)$ be such that for some $c_f > 0$ 
$$
\vert f(\zeta)\vert \leq 
c_f e^{\varphi(\zeta)+h(\zeta)}, \ 
\zeta\in\cal O.
$$ 

Then there exists a function $F\in H(\mathbb C^n \times \cal O)$ such that for $\zeta \in {\cal O}$ \ 
$F(\zeta,\zeta) = f(\zeta)$  and for some $C > 0$ and $N \geq 0$
$$
\vert F(z,\zeta)\vert\leq C (1 + \Vert (z, \zeta)\Vert)^N \left(1 + \frac {1} {\Delta_{\cal O}(\zeta)}\right)^N  e^{\varphi (z) + h(\zeta)},\ \ z\in\mathbb C^n,\ \ \zeta\in\cal O.
$$
\end{lemma}

{\bf Proof}. Let $\mu\in C^{\infty}[0,\infty)$ be such that for $t \geq 0$ \ $0\leq \mu(t)\leq 1$ and $\vert \mu'(t) \vert \leq 4$, 
$\mu (t)=1$ for $t\in [0,\frac {1}{3}], \mu(t)=0$ for $t\geq 1$. 
For ${\mathbb C} \times {\cal O}$ let
$
H_1 (z_1, \zeta) = \mu ((1+\Vert\zeta\Vert)^{\nu}\vert z_1\vert).
$
Note that $H_1 (z_1,\zeta)=0$ out of 
$$
\Omega_1 =
\{(z_1,\zeta)\in\mathbb C^{n+1}: (1+\Vert\zeta\Vert)^{\nu} 
\vert z_1 \vert < 1 \}.
$$ 
Let 
$$
W_1 = \{(z_1 ,\zeta) \in {\mathbb C}^{n+1}:
\frac {1}{3} < (1+\Vert \zeta \Vert)^{\nu} \vert z_1 \vert) < 1 \}.
$$
For $(z_1,\zeta)\in\mathbb C\times \cal O$ 
$$
\frac {\partial H_1} {\partial
\overline{\zeta}_j}(z_1 ,\zeta) = 
\frac {1}{2} 
\mu'((1 + \Vert \zeta \Vert)^{\nu} \vert
z_1 \vert) \nu \vert z_1 \vert (1 + \Vert \zeta \Vert)^{\nu -
1} \frac {\zeta_j} {\Vert \zeta \Vert},\ \ j=1,\ldots,n, 
$$
$$
\frac {\partial H_{1}} {\partial \overline{z}_1}(z_1 ,\zeta) = \frac {1} {2}
\mu'((1 + \Vert\zeta \Vert)^{\nu} \vert z_1 \vert) (1 +
\Vert\zeta \Vert)^{\nu} \frac {z_1} {\vert z_1 \vert}.
$$ 
Obviously, $\overline \partial H_1
(z_1 ,\zeta) = 0$ out of $W_1$. For $(z_1,\zeta)\in W_1$ 
$$\Vert \overline{\partial} H_1 (z_1,\zeta)\Vert ^{2} = 
\sum_{j=1}^{n}\left\vert \frac {\partial H_1(z_1, \zeta)} {\partial \zeta_j}\right\vert^2 + \left\vert\frac{\partial H_1 (z_1, \zeta)}{\partial\overline{z}_1}\right\vert^2 =
$$
$$
= \frac {1} {4}(\mu' ((1 + \Vert\zeta \Vert)^{\nu} \vert z_1 \vert)(1 +\Vert\zeta \Vert)^{\nu})^2 
\left(\frac {\nu^{2} \vert z_1\vert ^2} {(1 + \Vert \zeta \Vert)^2} + 1\right)\leq
$$
$$
\leq 4 (1 + \Vert \zeta \Vert)^{2\nu} \left(\frac {\nu^2} {(1 +\Vert\zeta\Vert)^{2 + 2 \nu}} +1\right) 
\leq 4(\nu^2 + 1) (1 +\Vert\zeta\Vert)^{2\nu}.
$$

Thus,  in $\mathbb C \times \cal O$ 
$$
\Vert
\overline{\partial} H_1 (z_1,\zeta)\Vert^2 \leq 4 (\nu^2
+1)(1+\Vert\zeta\Vert)^{2\nu}.
$$ 

Choose a function $v_1\in C^{\infty}(\mathbb C \times\cal O)$ with suitable bound so that 
$
F_1 (z_1, \zeta) = f(\zeta) H_1 (z_1, \zeta) - z_1 v_1(z_1, \zeta)
$  
is a holomorphic function in $\mathbb C \times \cal O$.
Obviously, $F_1(0,\zeta) = f(\zeta), \zeta \in \cal O$. 
The function $v_1\in C^{\infty}({\mathbb C} \times\cal O)$ will be found as a solution of the equation 
$$
\overline{\partial}v_1 (z_1, \zeta) = \frac{f(\zeta)\overline{\partial} H_1 (z_1, \zeta)} {z_1}, \ (z_1, \zeta) \in {\mathbb C} \times \cal O.
 \eqno (17)
$$
Denote by $g_1(z_1, \zeta)$ the form to the right of (17). 
It is easy to check that $\overline{\partial} g_1(z_1,\zeta)=0$  in $\mathbb C\times \cal O$. 
It is clear that
$g_1(z_1, \zeta)$ is a zero form of a type $(0,1)$ out of $W_1$. 
If
$(z_1, \zeta)\in W_1$ then 
$$
\Vert g_1 (z_1, \zeta)\Vert^2 =
\frac {\vert f(\zeta)\vert^2} {\vert z_1\vert^2} \Vert
\overline{\partial} H_1 (z_1, \zeta) \Vert^2 \leq c_f^2
e^{2(\varphi(\zeta) + h(\zeta))} 36 (\nu^2 + 1)(1+ \Vert \zeta
\Vert)^{4\nu}.
$$ 
Putting $B_1 =36 (\nu^2 + 1)c^2_f$ we have in $\mathbb C\times \cal O$ 
$$
\Vert g_1 (z_1, \zeta)\Vert ^2\leq B_1 e^{2(\varphi(\zeta) + h(\zeta))} (1 +\Vert\zeta\Vert)^{4\nu}.
$$

Recall that for $(z_1,\zeta)\in\mathbb C^{n+1}$ \ 
$\varphi_1(z_1, \zeta) = \varphi (\zeta_1 + z_1, \zeta_2,
\ldots, \zeta_n)$ 
and $\vert \varphi_1
(z_1,\zeta) - \varphi(\zeta)\vert\leq c_{\varphi}$ if $(z_1, \zeta) \in \Omega_1$. 

Further
$$
\int_{\mathbb C\times\cal O} 
\frac {\Vert g_1 (z_1, \zeta)\Vert^2 e^{-2(\varphi_1 (z_1, \zeta) + h(\zeta))}} 
{(1 + \Vert (z_1, \zeta)\Vert)^{2n+1+2\nu}} d\lambda_{n+1}(z_1, \zeta)= 
$$
$$
=\int_{W_1} \frac {\Vert g_1 (z_1, \zeta)\Vert^2 e^{-2(\varphi_1(z_1, \zeta) + h(\zeta))} } {(1 + \Vert (z_1,\zeta)\Vert)^{2n+1+2\nu}} d\lambda_{n+1}(z_1, \zeta)\leq
$$
$$
\leq B_1\int_{W_1} \frac {(1+\Vert\zeta\Vert)^{2\nu} e^{2(\varphi(\zeta) - \varphi_1(z_1,\zeta))}} {(1 + \Vert (z_1,\zeta)\Vert)^{2n+1}} d\lambda_{n+1}(z_1, \zeta)\leq 
$$
$$
\leq B_1 e^{2c_{\varphi}}\int_{W_1} \frac {(1+\Vert\zeta\Vert)^{2\nu}} {(1 + \Vert (z_1,\zeta)\Vert)^{2n+1}} d\lambda_{n+1}(z_1,\zeta)
\leq
$$
$$
\leq B_1 e^{2c_{\varphi}}\int_{\mathbb C^n}\frac {(1+\Vert\zeta\Vert)^{2\nu}} {(1+\Vert\zeta\Vert)^{2n+1}}
\left(\int_{\vert z_1\vert < \frac {1} {(1+\Vert\zeta\Vert)^{\nu}}}d\lambda_1 (z_1)\right) d\lambda_{n}(\zeta) =
$$
$$
= \pi B_1 e^{2C_{\varphi}} \int_{\mathbb C^n} \frac {d\lambda_n (\zeta)} {(1+\Vert\zeta\Vert)^{2n+1}} < \infty .
$$

Theorem $2.2.1'$ in [12] let us to find a solution 
$v_1\in C^{\infty}(\mathbb C\times \cal O)$ of the equation (17) such that
$$
2\int_{\mathbb C\times\cal O} \frac{\vert v_1(z_1, \zeta)\vert^2 
e^{-2(\varphi_1(z_1,\zeta) + h(\zeta))}} {(1+\Vert (z_1,\zeta)\Vert)^{2n+2\nu+1} (1+\Vert (z_1,\zeta)\Vert^2)^2}d\lambda_{n+1}(z_1,\zeta)\leq
$$
$$
\leq \int_{\mathbb C\times\cal O}\frac {\Vert g_1 (z_1,\zeta)\Vert^2 e^{-2(\varphi_1(z_1,\zeta) + h(\zeta))}} {(1+\Vert (z_1,\zeta)\Vert)^{2n+2\nu+1}}d\lambda_{n+1}(z_1,\zeta). 
$$
From this
$$
\int_{\mathbb C\times\cal O}
\frac {\vert F_1(z_1,\zeta)\Vert^2 e^{-2(\varphi_1(z_1,\zeta) + h(\zeta))}} {(1+\Vert (z_1,\zeta)\Vert)^{2n+2\nu+1} (1+ \Vert (z_1,\zeta)\Vert^2)^3} d\lambda_{n+1}(z_1,\zeta)\leq 
$$
$$
\leq 2\int_{\mathbb C\times\cal O}\frac {\vert f(\zeta)\vert^2 \vert H_1 (z_1,\zeta)\vert^2 e^{-2(\varphi_1(z_1,\zeta) + h(\zeta))}} {(1+\Vert (z_1,\zeta)\Vert)^{2n+2\nu+1} (1+ \Vert (z_1,\zeta)\Vert^2)^3} d\lambda_{n+1}(z_1,\zeta)+
$$
$$+2 \int_{\mathbb C\times\cal O}\frac {\vert z_1\vert^2 \vert v_1(z_1,\zeta)\vert^2 e^{-2(\varphi_1(z_1,\zeta) + h(\zeta))}}{(1+\Vert (z_1,\zeta)\Vert)^{2n+2\nu+1} (1+ \Vert(z_1,\zeta)\Vert^2)^3} d\lambda_{n+1}(z_1,\zeta)\leq $$
$$
\leq 2 c^2_{f} \int_{\Omega_1}\frac {e^{2(\varphi(\zeta) - \varphi_1(z_1,\zeta)))}} {(1+\Vert (z_1,\zeta)\Vert)^{2n+2\nu+1} (1+ \Vert(z_1,\zeta)\Vert^2)^3} d\lambda_{n+1}(z_1,\zeta)+  
$$
$$
+ 2 \int_{\mathbb C\times\cal O}\frac {\vert v_1 (z_1,\zeta)\vert^2 e^{-2(\varphi_1(z_1,\zeta) + h(\zeta))}} {(1+\Vert (z_1,\zeta)\Vert)^{2n+2\nu+1} (1+ \Vert(z_1,\zeta)\Vert^2)^2} d\lambda_{n+1}(z_1,\zeta)\leq
$$
$$
\leq 2 c^2_{f} e^{2c_{\varphi}}\int_{\mathbb C^{n+1}}\frac {d\lambda_{n+1}(z_1,\zeta)} {(1+ \Vert(z_1,\zeta)\Vert)^{2n+1}} +  
$$
$$
+ \int_{\mathbb C\times\cal O}\frac {\Vert g_1 (z_1,\zeta)\Vert^2 e^{-2(\varphi_1(z_1,\zeta) + h(\zeta))}} {(1+\Vert (z_1,\zeta)\Vert)^{2n+2\nu+1} (1+ \Vert(z_1,\zeta)\Vert^2)^3} d\lambda_{n+1}(z_1,\zeta) = B_2 < \infty .
$$

Now we pass to uniform estimates on $F_1$. 
Let 
$(z_1,\zeta)\in \mathbb C \times \cal O $ and $R=\min \left(\frac {\Delta_{\cal O}(\zeta)} {4}, \frac {1}
{2(1+\sqrt2\Vert(z_1,\zeta)\Vert)^{\nu}}\right)$. 
Since $\vert F_1(z_1,\zeta)\vert^2$ is a plurisubharmonic function in 
$\mathbb C\times \cal O$ then 
$$
\vert F_1 (z_1,\zeta)\vert^2\leq \frac {1} {\nu_{n+1}(R)}\int_{B_R(z_1,\zeta)} \vert F_1(t_1,w) \vert^2 d\lambda_{n+1}(t_1,w).
$$
From this and previous integral estimate on $F_1$ we obtain
$$
\vert F_1 (z_1,\zeta)\vert^2\leq \frac {B_2} {\nu_{n+1} (1)}\left(\frac {4} {\Delta_{\cal O}(\zeta)} +2(1+\sqrt 2\Vert (z_1,\zeta)\Vert^{\nu})\right)^{2(n+1)}\cdot 
$$
$$
\cdot(2+\Vert
(z_1,\zeta)\Vert)^{2n+2\nu+7}e^{\displaystyle\sup_{(t_1,w)\in B_R((z_1,\zeta),R)}2(\varphi_1(t_1,w)+h(w))}.
$$
Since here 
$\Vert w-\zeta\Vert \leq \min (1, \frac {\Delta_{\cal O}(\zeta)} {4})$ 
then 
$\vert h(w) - h(\zeta)\vert\leq c_h$. 
Also we have for 
$(t_1,w)\in B_R(z_1,\zeta)$  
$$
\Vert (t_1+w_1,w_2,\ldots,w_n) - (z_1+\zeta_1,\zeta_1,\ldots,\zeta_n)\Vert 
\le 
\vert t_1 - z_1 \vert + \Vert w - \zeta \Vert \leq 
$$ 
$$
\leq 2R \le 
\frac {1} {(1 + \sqrt 2 \Vert(z_1,\zeta)\Vert)^{\nu}}\leq \frac {1}
{(1+\Vert(z_1+\zeta_1,\zeta_2,\ldots,\zeta_n)\Vert)^{\nu}} \ .
$$
Hence for 
$(t_1,w)\in B_R(z_1,\zeta)$
$$
\vert \varphi_1 (t_1,w) - \varphi_1 (z_1,\zeta)\vert =
\vert \varphi(t_1+w_1,w_2,\ldots,w_n) - \varphi
(z_1+\zeta_1,\zeta_2,\ldots,\zeta_n)\vert\leq c_{\varphi}. 
$$
Using this inequalities we find a constant $C_1>0$ such that for $(z_1,\zeta)\in\mathbb C\times\cal O$
$$
\vert F_1(z_1,\zeta)\vert^2 \leq C_1 (1 + \Vert(z_1,\zeta)\Vert)^{n+(n+2)\nu+4} \left(1+\frac {1} {\Delta_{\cal O}(\zeta)}  \right)^{n+1}e^{\varphi_1(z_1,\zeta)+h(\zeta)}.
$$
Successively applying lemma 9 $(n-1)$ times we will construct 
a function $F_n\in H(\mathbb C^n \times \cal O)$  such that 
$F_n(0, \ldots, 0, \zeta)=f(\zeta)$ for $\zeta \in \cal O$ and for some $c_n>0$ and $N \geq 0$ 
$$
\vert
F_n(z, \zeta)\vert\leq
C_n(1+\Vert(z,\zeta)\Vert)^N \left(1 + \frac{1}
{\Delta_{\cal O}(\zeta)}\right)^N
e^{\varphi_n(z,\zeta) + h(\zeta)}, \ (z, \zeta) \in \mathbb C^n \times \cal O.
$$ 
For $(z,\zeta) \in {\mathbb C}^n \times {\cal O}$ 
put $F(z, \zeta)= F_n(z -\zeta, \zeta)$. Then $F\in H(\mathbb C^n \times \cal O)$, $F(\zeta,\zeta) = f(\zeta)$ 
for $\zeta \in {\cal O}$ and  
$$
\vert F(z, \zeta)\vert 
\leq C (1 + \Vert (z, \zeta)\Vert)^N \left(1 + \frac {1} {\Delta_{\cal O}(\zeta)}\right)^N  e^{\varphi (z) + h(\zeta)}, \ (z, \zeta) \in \mathbb C^n \times \cal O,
$$
where
$C=(1 + \sqrt 2)^N C_n$.

Lemma 10 is proved.

{\bf Remark 3}. The definition of functions $H_k$ in lemmas 9, 10 comes from [14].

\begin{lemma} 
Let $\cal O$ be a domain of holomorphy in $\mathbb C^n$. Let  $\varphi\in psh (\mathbb C^n)$ be such that for some
$c_{\varphi}>0$ and $\nu>0$ 
$
\vert\varphi(z)-\varphi(t)\vert\leq c_{\varphi}
$
if $\Vert z-t\Vert\leq\frac{1} {(1+\Vert t
\Vert)^{\nu}}$, 
$h\in psh(\cal O)$ be such that for some $c_h>0$
$
\vert h(z) - h(t) \vert \leq c_h
$
if for 
$z, t \in \cal O$ \ $\Vert z - t \Vert\leq \min
\left(1,\frac{\Delta_{\cal O}(t)} {4}\right)$. 

Let a function $S \in H({\mathbb C}^n \times {\cal O})$ satisfies to inequality
$$
\vert S(z, \zeta)\vert \leq 
e^{\varphi(z) + h(\zeta)}, \ z\in {\mathbb C}^n, 
\zeta \in {\cal O},
$$ 
and $S(\xi, \xi) = 0$ for $\xi \in {\cal O}$.

Then there exist functions $S_1, \ldots , S_n \in H({\mathbb C}^n \times {\cal O})$, numbers $C > 0$ and $N \ge 0$ such that for 
$(z, \xi)\in {\mathbb C}^n \times {\cal O}$:

a) $ S(z, \xi) = \displaystyle\sum_{j=1}^n S_j (z, \xi)(z_j - \xi_j)$;
 
b) for each $j=1, \ldots , n$
$$
\vert S_j(z, \xi) \vert \leq C (1 + \Vert (z, \zeta)\Vert)^N \left(1 + \frac {1} {\Delta_{\cal O}(\zeta)}\right)^N  e^{\varphi (z) + h(\zeta)}.
$$
\end{lemma}

{\bf Proof}. 
Let $L(z, \xi) = S(z+\zeta, \zeta)$, \ $z \in {\mathbb C^n}, \zeta \in \cal O.$ 
Then
$$
\vert L(z, \zeta)\vert \leq e^{\varphi_n(z,\zeta) + h(\zeta)}, \  (z, \zeta)\in {\mathbb C}^n \times {\cal O}, \eqno (18)
$$ 
and 
$L(0, \zeta) = S(\zeta, \zeta) = 0$ for $\zeta \in {\cal O}$.

For $z_1, \ldots , z_n \in {\mathbb C}, z \in {\mathbb C}^n, \zeta \in {\cal O}$ let
$L_1 (z_1,\zeta) = L(z_1,0,\ldots,0,\zeta)$,  
$L_2 (z_1,z_2,\zeta) = L(z_1,z_2,0,\ldots,0,\zeta)$, $ \ldots $ , 
$L_n(z, \zeta) = L(z,\zeta)$.
In view of (18) for  $k=1, \ldots , n$
$$
\vert L_k(z, \zeta) \vert \le e^{\varphi_k(z,\zeta) + h(\zeta)}, \  (z, \zeta)\in {\mathbb C}^k \times {\cal O}. \eqno (19)
$$

Since $L_1(0,\zeta)=0  \ \forall \zeta \in {\cal O}$ then 
$$
\psi_1^{(1)}(z_1,\zeta) = \frac{L_1(z_1, \zeta)}{z_1}
$$
is a holomorphic function in ${\mathbb C} \times {\cal O}$. 
Let us  estimate its growth. 
For  $\zeta\in {\cal O}$,
$\vert z_1 \vert \geq \frac {1}{(1 + \Vert \zeta \Vert)^{\nu}}$ 
$$
\vert\psi_1^{(1)}(z_1, \zeta)\vert \leq (1 + \Vert \zeta \Vert)^{\nu} e^{\varphi_1(z,\zeta) + h(\zeta)}. 
$$
For $\zeta \in {\cal O}$, $\vert z_1 \vert < \frac {1}{(1 + \Vert \zeta \Vert)^{\nu}}$ 
$$
\vert\psi_1^{(1)}(z_1, \zeta)\vert 
\leq 
\max \limits_{\vert t_1\vert = \frac {1}{(1 + \Vert \zeta \Vert)^{\nu}}} 
\left \vert \frac{L_1(t_1,\zeta)}{t_1} \right \vert
\le
e^{2 c_{\varphi}} (1 + \Vert \zeta \Vert)^{\nu} e^{\varphi_1(z_1, \zeta) + h(\zeta)}.
$$
Thus, putting 
$A_1 = e^{2 c_{\varphi}}, m_1=\nu$ for $z_1 \in {\mathbb C}, \zeta \in {\cal O}$
$$
\vert \psi_1^{(1)}(z_1,\zeta)\vert\leq A_1 (1 + \Vert (z_1, \zeta) \Vert)^{m_1} 
\left(1 + \frac {1} {\Delta_{\cal O}(\zeta)}\right)^{m_1} e^{\varphi_1(z_1, \zeta) + h(\zeta)}.
$$
Note also that in  ${\mathbb C} \times {\cal O}$ \  $L_1(z_1,\zeta) = \psi_1^{(1)}(z_1, \zeta)z_1$. 

Let for $k=2,\ldots,n$ there are functions $\psi_j^{(k-1)} \in H({\mathbb C}^{k-1} \times {\cal O})$ ($j=1, \ldots, k-1$) such that
$$
L_{k-1}(z_1, \ldots, z_{k-1}, \zeta) = \displaystyle \sum_{j=1}^{k-1} \psi_j^{(k-1)}(z_1,\ldots,z_{k-1}, \zeta)z_j
$$
and numbers $A_{k-1}>0$ and $m_{k-1}\in {\mathbb N}$ such that for all $z' \in {\mathbb C^{k-1}}, \zeta\in {\cal O}, j= 1, \ldots , k-1$
$$
\vert \psi_j^{(k-1)}(z',\zeta)\vert\leq A_{k-1} (1 + \Vert (z', \zeta) \Vert)^{m_{k-1}} 
\left(1 + \frac {1} {\Delta_{\cal O}(\zeta)}\right)^{m_{k-1}} e^{\varphi_{k-1}(z', \zeta) + h(\zeta)}.
$$
According to lemma 9 for each $j=1,\ldots,k-1 $ there exists a function $\psi_j^{(k)}$ holomorphic in ${\mathbb C^k} \times {\cal O}$ such that 
$$
\psi_j^{(k)}(z_1, \ldots, z_{k-1}, 0, \zeta) = \psi_j^{(k-1)}(z_1, \ldots, z_{k-1}, \zeta), \ z_1, \ldots, z_{k-1} \in {\mathbb C}, 
\zeta \in {\cal O}, 
$$
and numbers $B_k>0$ and $\tilde m_k \geq 0$ such that for $z = (z', z_k) \in \mathbb C^k, \zeta \in \cal O$ and for all $j=1,\ldots,k-1 $
$$
\vert \psi_j^{(k)}(z, \zeta) \vert \leq B_k(1+\Vert(z, \zeta)\Vert)^{\tilde m_k} \left(1 + \frac{1}
{\Delta_{\cal O}(\zeta)}\right)^{\tilde m_k}
e^{\varphi_k(z, \zeta) + h(\zeta)}.  \eqno (20)
$$ 
Put
$$
Y_k(z, \zeta)=L_k(z, \zeta) -\psi_1^{(k)}(z, \zeta) z_1 - \ldots - \psi_{k-1}^{(k)}(z, \zeta)z_{k-1}, \ z \in\mathbb C^k,  \zeta \in \cal O.
$$
Using (19) and (20) we find a constant $C_k > 0$ such that 
$$
\vert Y_k(z, \zeta) \vert \leq C_k (1+\Vert(z, \zeta)\Vert)^{\tilde m_k+1} \left(1 + \frac{1}
{\Delta_{\cal O}(\zeta)}\right)^{\tilde m_k}
e^{\varphi_k(z, \zeta) + h(\zeta)}. 
$$ 
Also note that for $z_1, \ldots , z_{k-1} \in \mathbb C, \zeta \in  \cal O$ \ 
$
Y_k(z_1, \ldots,z_{k-1}, 0, \zeta)= 0.
$
Hence 
$$
\psi_k^{(k)}(z, \zeta)=\frac{Y_k(z, \zeta)} {z_k} 
$$
is a holomorphic function in $\mathbb C^k\times \cal O.$

Let us estimate a growth of  $\psi_k^{(k)}$. 
Put
$R_k=\frac {1} {2(1+\sqrt2\Vert(z', \zeta)\Vert)^{\nu}}$. 
For $z^{'}\in\mathbb C^{k-1}$, $\zeta\in\cal O$, $\vert z_k \vert \geq R_k$ 
$$
\vert\psi_k^{(k)}(z, \zeta)\vert < 2^{\nu+1} C_k (1+\Vert(z, \zeta)\Vert)^{\tilde m_k + \nu +1} \left(1 + \frac{1}
{\Delta_{\cal O}(\zeta)}\right)^{\tilde m_k}
e^{\varphi_k(z, \zeta) + h(\zeta)}. 
$$
For $z' \in\mathbb c^{k-1}$, $\zeta\in\cal O$, $\vert z_k \vert < R_k$, 
$$
\vert \psi_k^{(k)}(z', z_k, \zeta)\vert 
\leq 
\max \limits_{\vert t_k \vert = R_k} 
\left \vert \frac{Y_k(z', t_k, \zeta)}{t_k} \right \vert \le
$$
$$
\le 2^{\nu+ \tilde m_k +2} \left(1+\frac{1} {\Delta_{\cal O}(\zeta)}  \right)^{\widetilde{m_k}}
e^{2 c_{\varphi}} (1 + \Vert (z,  \zeta) \Vert)^{\tilde m_k + \nu +1} e^{\varphi_k(z, \zeta) + h(\zeta)}.
$$
Put $m_k=\widetilde{m}_k + \nu + 1$. From this and (20) it follows that there exists a constant 
$A_k >0$ such that for all $j=1,\ldots, k $ 
$$
\vert \psi_j^{(k)}(z,\zeta) \vert \leq A_k (1 + \Vert (z, \zeta) \Vert)^{m_k} 
\left(1 + \frac {1} {\Delta_{\cal O}(\zeta)}\right)^{m_k} e^{\varphi_k(z, \zeta) + h(\zeta)}. 
$$
Note that 
$
L_k(z, \zeta) = \psi_1^{(1)}(z, \zeta)z_1+\ldots+\psi_k^{(k)}(z, \zeta)z_k, \ z \in \mathbb C^k , \zeta \in \cal O. 
$
So if $k=n$ then
$$
L_n(z, \zeta)=L(z, \zeta)=\psi_1^{(n)}(z, \zeta)z_1+\ldots+\psi_n^{(n)}(z, \zeta)z_n, \ z\in\mathbb C^n, \zeta\in \cal O.
$$
Besides that for each $j=1, \ldots, n$
$$
\vert \psi_j^{(n)}(z,\zeta) \vert\leq A_n (1 + \Vert (z, \zeta) \Vert)^{m_n} 
\left(1 + \frac {1} {\Delta_{\cal O}(\zeta)}\right)^{m_n} e^{\varphi_n(z, \zeta) + h(\zeta)}. \eqno (21)
$$
Thus, we have 
$$
S(z,\zeta)=\psi_1^{(n)}(z-\zeta,\zeta)(z_1-\zeta_1)+\ldots+\psi_n^{(n)}(z-\zeta,\zeta)(z_n-\zeta_n),z\in\mathbb C^n, \zeta\in \cal O. 
$$
For $z\in\mathbb C^n, \zeta\in \cal O$ let $S_j(z, \zeta)=\psi_j^{(n)}(z-\zeta, \zeta)$. 
Then
$$
S(z, \zeta)=\sum_{j=1}^n S_j(z,\zeta)(z_j-\zeta_j),\ \ z\in\mathbb C^n,\ \ \zeta\in \cal O.
$$
Putting $N=m_n$, $C=A_n 2^{m_n}$ from  (21)  we have for each $j=1,\ldots,n$
$$
\vert S_j(z,\zeta)\vert\leq 
C(1 + \Vert (z, \zeta)\Vert)^N \left(1 + \frac {1} {\Delta_{\cal O}(\zeta)}\right)^N  e^{\varphi (z) + h(\zeta)}, \ z\in\mathbb C^n, \zeta\in \cal O.
$$

Lemma 11 is proved.

{\bf Remark 4}. The idea of proving lemma 11 comes from [4].

{\bf Proof of Theorem 3}. By lemma 8 the linear mapping $L: S \in G_M^*(U) \rightarrow \hat S$ acts from $G_M^*(U)$ into $H_{b, M}(T_C)$. 

Before we show that $L$ is continuous note that the topology of $G_M^*(U)$ can be described as follows.
Let $W_k = \{ f \in G_M(U):  \ p_k(f) \leq 1 \},  \ k \in {\mathbb N}$. Let 
$W_k^0 = \{ F \in G_M'(U): \vert (F, f) \vert \leq 1,  \ \forall f \in W_k \}$ be a polar  of $W_k$ in $G_M'(U)$. Let 
$E_k = \displaystyle \bigcup_{\alpha >0} (\alpha W_k^0) $
be a vector subspace in $G_M'(U)$ generated by polar $W_k^0, k =1,2,\ldots$. 
Define a topology in $E_k$ with the help of the norm
$$q_k(F) = \displaystyle \sup_{f \in W_k} \vert (F, f) \vert,  \ F\in E_k. $$
Note that  
$G_M'(U) = \bigcup_{k=1}^\infty E_k.$ 
Define in $G_M'(U)$ the topology $\lambda$ of an inductive limit of spaces $E_k$. Since $G_M(U)$ is a reflexive space 
then the strong topology in $G_M'(U)$ coincides with the topology $\lambda$ [15, chapter 8]. 

Now let $S \in E_m, \ m \in{\mathbb N}$. 
Then $$\vert (S, f) \vert \leq q_m(S), \  f \in W_m.
$$ 
From this it follows that 
$$
\vert (S, f) \vert \leq q_m(S) p_m(f), \ f \in G_M(U). 
$$
Putting here $f(\xi) = \exp(i<\xi, z>)$  with $z$ in $T_C$ and using (14) we obtain  
$$
\vert \hat S(z) \vert  \leq q_m(S) A e^{\omega_{m + [r] + 1} (\Vert z \Vert)} \left(1 + \frac {1}{\Delta_C(y)}\right)^{3m} \ ,  \eqno (22)
$$
where a constant $A>0$ does not depend on $z \in T_C$. Let $N(m)= \max (m + [r] + 1, 3m)$. Then from  (22) it
follows that 
$$
{\Vert \hat S\Vert}_{N(m)} \leq  A q_m(S),  \ S \in E_m  \ (m =1, 2, \ldots).
$$ 
Thus, $L$ is continuous.

While proving that $L$ is bijective we follow a scheme from [4].

First  show that $L$ is surjective. 
Let $F \in H_{M}(T_C)$. This means that $F$ is holomorphic in $H(T_C)$ and for some $c>0, m \in {\mathbb N}$ 
$$
\vert F(z) \vert \leq c e^{\omega_m (\Vert z \Vert)} \left (1 + \frac {1} {\Delta_C (y)}\right)^m. 
$$
Since for some $r > e $ \  $\vert b(y) \vert \le r \Vert y \Vert$ ($y \in {\overline C}$)  then
$$
\vert F(z) \vert \leq c e^{b(y) + \omega_m (\Vert z \Vert) + r \Vert z \Vert} \left (1 + \frac {1} {\Delta_C (y)}\right)^m, \ z \in T_C.
$$
Using (1) and putting $c_1 = c e^{([r]+1)Q}$, $k=m+[r]+1$ we have 
$$
\vert F(z) \vert \leq c_1 e^{\omega_k(\Vert z \Vert)} e^{b(Im z)} \left (1 + \frac {1} {\Delta_{T_C} (z)}\right)^m, \ z \in T_C. 
$$

Note that functions $h(\zeta) = b(Im \zeta) + m \ln (1 + \frac {1} {\Delta_{T_C} (\zeta)})$ $ (z\in\mathbb C^n, \zeta \in T_C$) and 
$\varphi (z) = \omega_k (\Vert z \Vert)$  satisfy to conditions of lemma 10. Then there exists a function $\Phi \in H({\mathbb C^n} \times T_C)$ such that $\Phi(\zeta, \zeta) = F(\zeta)$ for $\zeta \in T_C$ 
and numbers $c_2 > 0$ and $N \geq 0$ such that for $z \in\mathbb C^n,   \zeta \in T_C$
$$
\vert \Phi(z,\zeta)\vert\leq c_2 (1 + \Vert (z, \zeta)\Vert)^N \left(1 + \frac {1} {\Delta_{T_C} (\zeta)} \right)^{N} e^{b(Im \zeta)}
e^{\omega_k (\Vert z \Vert)}. \eqno (23)
$$

Since $\Phi(z, \zeta)$ is an entire function of $z$ then it can be expanded in a power series:
$$
\Phi(z, \zeta) = \displaystyle \sum_{\vert\alpha\vert\geq 0} C_\alpha (\zeta) z^\alpha,  \ \zeta\in T_C,  \ z\in {\mathbb C^n}.
$$
By the Cauchy formula
$$
C_\alpha (\zeta) = \frac{1}{(2\pi i)^n} \int_{\vert z_1\vert =R}\ldots\int_{\vert z_n\vert =R} \frac {\Phi (z, \zeta)} {z_1^{\alpha_1 +1} \ldots z_n^{\alpha_n +1}} \ dz_1\ldots dz_n,$$
where $\alpha\in {\mathbb Z_+^n, \ \ R>0}$ is arbitrary. 
From this it follows that $C_\alpha \in H(T_C)$. 
Using (23)  we have for $\zeta \in T_C$
$$
\vert C_\alpha(\zeta) \vert \leq 
\frac 
{c_2 (1 + \sqrt{n} R)^N(1 + \Vert \zeta \Vert)^N e^{b(Im \zeta) + \omega_k(\sqrt{n} R)}
\left(1 + \frac {1} {\Delta_{T_C} (\zeta)} \right)^{N}} {R^{\vert \alpha\vert}} \ .
$$
Using lemma 1 one can find a constant $c_3>0$ such that for each
$  R>0$
$$
\vert C_\alpha (\zeta)\vert \leq 
c_3 \frac 
{e^{\omega_{k+1}(\sqrt{n} R)}} {R^{\vert\alpha\vert}} (1 + \Vert \zeta \Vert)^N
e^{b(Im \zeta)} \left(1 + \frac {1} {\Delta_{T_C} (\zeta)} \right)^{N}, \ \zeta \in T_C.
$$
Hence for $\zeta \in T_C$
$$
\vert C_\alpha (\zeta)\vert \leq c_3
\left(\frac{\sqrt{n}} {\varepsilon_{k+1}}\right)^{\vert\alpha\vert} \left(\inf_{R>0} \frac {e^{\omega_M(r)}} {r^{\vert\alpha\vert}}\right) 
e^{b(Im \zeta)} (1 + \Vert\zeta\Vert)^N \left(1 + \frac {1} {\Delta_{T_C} (\zeta)} \right)^{N}.
$$
Since [16]
$$
\displaystyle \inf_{R>0} \frac {e^{\omega_M(r)}} {r^k} = \frac{1} {M_k} \ ,  \ k = 0,1, \ldots,
$$
then for  $\alpha\in\mathbb Z_{+}^{n}$, $\zeta\in T_C$ we have
$$
\vert C_\alpha (\zeta)\vert \leq 
c_3 \left(\frac{\sqrt{n}} {\varepsilon_{k+1}}\right)^{\vert\alpha\vert} 
\frac{e^{b(Im \zeta)}} {M_{\vert\alpha\vert}} 
(1 + \Vert\zeta\Vert)^N \left(1 + \frac {1} {\Delta_{T_C} (\zeta)} \right)^{N}. \eqno (24)
$$
Thus, for each $\alpha\in {\mathbb Z_+^n}$ \ $C_\alpha \in V(T_C)$. 
By theorem 2 there exist functionals $S_\alpha \in S^*(U)$ such that $\hat {S_\alpha} = C_\alpha$. 

From (24) and properties of  $(LN^*)$-spaces it follows that  
$\{M_{\vert\alpha\vert} 
\left(\frac {\varepsilon_{k+1}} {\sqrt{n}} \right)^{\vert\alpha\vert} C_\alpha \}_{\alpha \in \mathbb Z_+^n} $ 
is a bounded set in $V(T_C)$. 
In view of the topological isomorphism of spaces $S^*(U)$ and $V(T_C)$ the set
${\cal A}= \{M_{\vert\alpha\vert} \left(\frac {\varepsilon_{k+1}} 
{\sqrt{n}} \right)^{\vert\alpha\vert} S_\alpha \}_{\alpha \in \mathbb Z_+^n}$
is bounded in $S^* (U)$. Hence it is  weakly bounded. 
By Schwartz theorem [1] there exist numbers $c_4 >0$ and $p\in\mathbb N$ such that
$$
\vert (F, f) \vert \leq c_4 {\Vert f \Vert}_{p, U}, \ F \in {\cal A}, \ f \in S(U).  
$$ 
Thus, for each
$\alpha \in {\mathbb Z_+^n}, f \in S(U)$ 
$$
\vert (S_{\alpha}, f) \vert \leq c_4 \left(\frac{\sqrt{n}}{\varepsilon_{k+1}}\right)^{\vert\alpha\vert} 
\frac{{\Vert  f \Vert}_{p, U}}{M_{\vert\alpha\vert}} \ . \eqno (25)
$$

Define a functional $T$ on $G_M(U)$ by the rule:
$$
(T, f) = \sum_{\vert\alpha\vert \geq  0} (S_\alpha, (-i)^{\vert\alpha\vert} D^\alpha f), \ f \in G_M(U).  \eqno (26)
$$
It is correctly defined. 
Indeed,
using (25)  we have for each $ f \in G_M(U), \alpha \in {\mathbb Z_+^n}, s \in {\mathbb N}$
$$
\vert (S_\alpha, D^\alpha f) \vert \leq 
c_4 \left(\frac {\sqrt{n}} {\varepsilon_{k+1}}\right)^{\vert\alpha\vert} \frac{1} {M_{\vert\alpha\vert}}\sup_{x \in U, \vert\beta\vert\leq p }\vert (D^{\alpha+\beta} f)(x)\vert (1 + \Vert x\Vert)^p\leq 
$$
$$
\leq c_4 \left(\frac {\sqrt{n}} {\varepsilon_{k+1}}\right)^{\vert\alpha\vert} \frac{1} {M_{\vert\alpha\vert}}\sup_{x\in U, \vert\beta\vert\leq p } \frac {p_s (f) \varepsilon_s^{\vert\alpha\vert +\vert\beta\vert} M_{\vert\alpha\vert +\vert\beta\vert}(1 + \Vert x\Vert)^p}  {(1 + \Vert x\Vert)^s} \ .
$$
Using the condition $i_2)$ for $s\geq p$ we have
$$
\vert (S_\alpha, D^\alpha f) \vert \leq c_4 
\left(\frac {\sqrt{n}}{\varepsilon_{k+1}}\right)^{\vert\alpha\vert} p_s (f) 
\sup_{\vert\beta\vert\leq p} \varepsilon_s^{\vert\alpha\vert + \vert\beta\vert} 
H_1 H_2^{\vert\alpha\vert + \vert\beta\vert} M_{\vert\beta\vert}\leq 
$$
$$
\leq c_4 \left(\frac {\sqrt{n}}{\varepsilon_{k+1}}\right)^{\vert\alpha\vert} H_1  M_p H_2^{\vert \alpha \vert} p_s (f) \varepsilon_s^{\vert\alpha\vert}
= c_4 H_1  M_p \left(\frac {\sqrt{n}\varepsilon_s H_2 }  {\varepsilon_{k+1}}\right)^{\vert\alpha\vert} p_s (f). 
$$
Now choose $s$ so that
$\tau_s=\frac {\sqrt{n}\varepsilon_s H_2 }  {\varepsilon_{k+1}}  <1$. 
Then for each $f \in G_M(U), \alpha \in {\mathbb Z_+^n}$ \ 
$$
\vert (S_\alpha, D^\alpha f) \vert \leq c_5 \tau_s^{\vert\alpha\vert} p_s (f), 
$$
where $c_5= c_4 H_1  M_p$.
From this it follows that the series to the right of (26) converges and that 
$$
\vert (T, f) \vert\leq \frac {c_4 }{(1 - \tau_s)^n} p_s (f),\ \ f \in G_M(U). 
$$
Hence the linear functional $T$ is correctly defined and continuous. 
Besides that $\hat T  = F.$ Indeed, for each $z \in T_C$
$$
\hat  T(z)= \sum_{\vert\alpha\vert \ge 0} (S_{\alpha}, (-i)^{\vert\alpha\vert} D^{\alpha} (e^{i<\zeta, z>})) 
=
\sum_{\vert\alpha\vert \ge 0}(S_{\alpha}, (-i)^{\vert\alpha\vert}(iz)^{\alpha}(e^{i<\zeta, z>}))
$$
$$
=
\sum_{\vert\alpha\vert \ge 0}z^{\alpha}(S_{\alpha},(e^{i<\zeta, z>}))=
\sum_{\vert\alpha\vert \ge 0} C_{\alpha}(z) z^{\alpha} = \Phi(z, z) = F(z).
$$
Thus, $L$ is surjective.

The mapping $L$ is injective. Indeed, let for $T \in G_M'(U)$ \ $\hat{T} \equiv~0$. We will show that $T$ is a zero functional. 
Since $T$ is a linear continuous functional there
exist numbers $m \in \mathbb N$ and $c_T>0$ such that
$$
\vert (T, f) \vert \leq c_T p_m (f), \ f \in G_M(U).  
$$
By lemma 7 there exist functionals $T_{\alpha}\in C_m'(U)$ ($\alpha\in\mathbb Z_+^n$) such that
$$
(T, f) = \displaystyle \sum_{\alpha \in \mathbb Z_+^n} (T_{\alpha}, D^{\alpha} f), \ f \in G_M(U),
$$ 
and 
$$
\vert (T_{\alpha}, g) \vert \leq \frac {c_T}  {\varepsilon_m^{\vert\alpha\vert} M_{\vert\alpha\vert}} \widetilde{p}_m (g), \ 
g \in C_m(U).  \eqno (27)
$$
From this we have for each $z \in T_C$
$$
\hat T(z) = \sum_{\alpha \in \mathbb Z_+^n} (T_{\alpha}, (iz)^{\alpha} e^{i<\zeta, z>}) = 
\sum_{\alpha \in \mathbb Z_+^n} i^{\vert\alpha\vert} (T_{\alpha}, e^{i<\zeta, z>}) z^{\alpha}.
$$
Let $V_{\alpha}(z) = i^{\vert\alpha\vert} (T_{\alpha}, e^{i<\xi, z>}) $. 
Obviously $V_{\alpha}\in H(T_C)$.
Using (27) and lemma 6 we obtain the estimate
$$
\vert V_{\alpha} (z) \vert \leq \frac {d_1} {\varepsilon_m^{\vert\alpha\vert} M_{\vert\alpha\vert}}(1 + \Vert z \Vert)^{2m}
\left(1 + \frac {1}{\Delta_C(y)}\right)^{3m} e^{b(y)}, \eqno (28)
$$
where $d_1 >0$ is some constant not depending on $z = x + iy  \in T_C$ and $\alpha\in\mathbb Z_+^{n}$. 
Consider now the function
$
S (u, z) = \displaystyle \sum_{\vert\alpha\vert\geq 0} V_{\alpha} (z) u^{\alpha}, \ z \in T_C, u \in \mathbb C^n.  
$
Using (28) we have
$$
\vert S(u, z)\vert  
\leq   \sum_{\vert \alpha\vert \geq 0} \frac {d_1 \Vert u\Vert^{\vert\alpha\vert}}  {\varepsilon_m^{\vert\alpha\vert} M_{\vert\alpha\vert}} (1 + \Vert z \Vert)^{2m} \left(1 + \frac {1}{\Delta_C(y)}\right)^{3m} e^{b(y)}=  
$$
$$
= d_1 \left(1 + \frac {1}{\Delta_C(y)}\right)^{3m} e^{b(y)}
\sum_{\vert \alpha\vert \geq 0}\frac {\Vert u\Vert^{\vert\alpha\vert}}  {\varepsilon_{m+1}^{\vert\alpha\vert} M_{\vert\alpha\vert}} \left(\frac {\varepsilon_{m+1}}  {\varepsilon_m}\right)^{\vert\alpha\vert} \leq
$$
$$
\leq d_1 e^{b(y)}\left(1 + \frac {1}{\Delta_C(y)}\right)^{3m}
\sup_{\alpha\in\mathbb Z_+^n} \frac {\Vert u\Vert^{\vert\alpha\vert}}  {\varepsilon_{m+1}^{\vert\alpha\vert} M_{\vert\alpha\vert}} \sum_{\vert\alpha\vert\geq 0} \left(\frac {\varepsilon_{m+1}}  {\varepsilon_m}\right)^{\vert\alpha\vert} =  
$$
$$
\leq d_1 e^{b(y)}\left(1 + \frac {1}{\Delta_C(y)}\right)^{3m}
e^{\omega_{m+1}(\Vert u\Vert)} \left(\frac {\varepsilon_m}  {\varepsilon_m - \varepsilon_{m+1}}\right)^n \ .
$$
Note that for each $z\in T_C$
$S(z, z) = \sum_{\vert\alpha\vert\geq 0} V_{\alpha}(z) z^{\alpha} = 0$. 
Then by lemma 11 there exist functions
$S_1,\ldots,S_n \in H(\mathbb C^n \times T_C)$ such that
$$
S(z, \zeta) = \sum_{j=1}^n S_j(z, \zeta) (z_j-\zeta_j), \ z \in \mathbb C^n, \zeta \in T_C,
$$ 
and numbers $d_2>0$ and $N \in \mathbb N$  such that for   $j = 1,2,\ldots,n$, $z \in \mathbb C^n, \zeta \in T_C$
$$
\vert S_j(z, \zeta)\vert \leq d_2 (1 + \Vert (z, \zeta)\Vert)^N \left(1 + \frac {1} {\Delta_{T_C}(\zeta)} \right)^{N}
e^{\omega_{m+1}(\Vert z\Vert) + b(Im \zeta)}. \eqno (29)
$$
Next, we expand  $S_j$ in powers of $z$:
$$
S_j(z, \zeta) = \sum_{\vert\alpha\vert\geq 0} S_{j, \alpha} (\zeta) z^{\alpha},\ \ z\in\mathbb C^n,\ \ \zeta\in T_C.  
$$
Using (29), (1) and Cauchy's inequality for coefficients of power series we have  
$$
\vert S_{j, \alpha}(\zeta) \vert 
\leq \inf_{R>0}\frac {\displaystyle\max_{\vert z_1\vert=R,\ldots,\vert z_n\vert=R} \vert S_j(z, \zeta)\vert}  {R^{\vert\alpha\vert}}\leq 
$$
$$
\leq \inf_{R>0} \frac {d_3e^{\omega_{m+2}(\sqrt{n}R)} 
e^{b(Im \zeta) + (N+3m) \ln (1 + \frac {1} {\Delta_{T_C}(\zeta)}) + N \ln (1+\Vert\zeta\Vert))}}  {R^{\vert\alpha\vert}} =  
$$
$$
= d_3 e^{b(Im\zeta)}e^{N \ln (1 + \frac {1} {\Delta_{T_C}(\zeta)}) +  N \ln (1+\Vert\zeta\Vert)} \inf_{R>0}\frac {e^{\omega (R)} \sqrt{n}^{\vert\alpha\vert}}  {(R\varepsilon_{m+2})^{\vert\alpha\vert}} = 
$$
$$
=d_3 e^{b(Im\zeta)}e^{N \ln (1 + \frac {1} {\Delta_{T_C}(\zeta)}) +  N \ln (1+\Vert\zeta\Vert)} \frac {1} {M_{\vert\alpha\vert}} 
\left(\frac {\sqrt{n}} {\varepsilon_{m+2}}\right)^{\vert\alpha\vert} \ ,   
$$
where $d_3>0$ is some constant.
Choose $k\in\mathbb N$ so that $\varepsilon_k\sqrt{n} < \varepsilon_{m+2}$. 
Then
$$
\vert S_{j,\alpha}(\zeta)\vert \leq \frac {A_1 e^{b(Im\zeta)}e^{N \ln (1 + \frac {1} {\Delta_{T_C}(\zeta)}) +  N \ln (1+\Vert\zeta\Vert)}}  {\varepsilon_k^{\vert\alpha\vert} M_{\vert\alpha\vert}}, \ \zeta\in T_C. \eqno (30)
$$
By theorem 2 there exist functionals $\psi_{j, \alpha} \in S^*(U)$ such that $\hat \psi_{j, \alpha} = S_{j, \alpha}$. 
By (30) it follows that the set 
$\{S_{j, \alpha}\varepsilon_k^{\vert\alpha\vert} M_{\vert\alpha\vert} \}_{\alpha \in \mathbb Z_+^n} $ 
is bounded in $V_b(T_C)$. 
But then the set 
$\Psi = \{\varepsilon_k^{\vert\alpha\vert} M_{\vert\alpha\vert}\psi_{j, \alpha} \}_{\alpha \in \mathbb Z_+^n, j=1,\ldots,n}$ 
is bounded in $S^*(U)$. Hence it is weakly bounded. 
By Schwartz's theorem [1, p. 93] there exist numbers $d_4>0$ and $p\in\mathbb N$ such that
$$
\vert (F, \varphi) \vert \leq d_4 \vert \Vert \varphi \Vert\vert_{p, U}, \ F \in \Psi, \ \varphi\in S(U).  
$$
Thus, for each $j=1,\ldots,n$ and $\alpha\in\mathbb Z_+^n$ 
$$
\vert (\Psi_{j, \alpha}, f) \vert\leq \frac{d_4} {\varepsilon_k^{\vert\alpha\vert} M_{\vert\alpha\vert}} {\Vert f\Vert}_{p, U}, \ 
f \in S(U). \eqno (31)
$$ 
For $j=1,\ldots,n$  and $\alpha \in \mathbb Z^n$ with at least one negative component 
let $\Psi_{j, \alpha}$ be a zero functional on $S(U)$ and $S_{j, \alpha} (z)=0,\ \ \forall z \in \mathbb C^n$. 
Then
$$
S(z, \zeta) = \sum_{j=1}^n S_j(z, \zeta) (z_j - \zeta_j) = 
\sum_{j=1}^n \sum_{\vert\alpha\vert\geq 0} S_{j, \alpha}(\zeta) z^{\alpha} (z_j - \zeta_j)=
$$ 
$$
=\sum_{j=1}^n \sum_{\vert\alpha\vert\geq 0} S_{j,\alpha}(\zeta)z_1^{\alpha_1}\ldots z_j^{\alpha_j +1} \ldots z^{\alpha_n}-\sum_{j=1}^n \sum_{\vert\alpha\vert\geq 0}S_{j,\alpha}(\zeta)z^{\alpha}\zeta_j=  $$ 
$$
=\sum_{j=1}^n \sum_{\vert\alpha\vert\geq 0} (S_{j,(\alpha_1,\ldots,\alpha_{j}-1,\ldots,\alpha_n)}(\zeta) - S_{j, \alpha}(\zeta) \zeta_j)z^{\alpha},\ \ z\in\mathbb C^n, \ \zeta\in T_C. 
$$ 
Hence
$$
V_{\alpha}(\zeta) = \sum_{j=1}^n (S_{j,(\alpha_1,\ldots,\alpha_{j-1},\ldots,\alpha_n)} (\zeta) - S_{j, \alpha}(\zeta)\zeta_j). \eqno (32)
$$ 
The expression to the right of (32) can be represented in the following form
$$
\sum_{j=1}^n (\hat \Psi _{j,(\alpha_1,\ldots,\alpha_{j-1},\ldots,\alpha_n)}(\zeta) + 
i (\Psi_{j, \alpha}, \frac {\partial} {\partial \xi_j}(e^{i<\xi, \zeta>})) ).
$$  
From this and theorem 2 it follows that 
$$
(T_{\alpha}, f) = (-i)^{\vert\alpha\vert} \sum_{j=1}^n (i(\Psi_{j, \alpha}, (\frac {\partial} {\partial\xi_j} f)) + (\Psi_{j, (\alpha_1, \ldots, \alpha_{j-1},\ldots, \alpha_n)}, f )), \ f \in S(U).
$$
So, for $f \in G_M(U)$
$$
(T, f) = \sum_{\vert\alpha\vert\geq 0} (T_{\alpha}, D^{\alpha} f) = 
$$
$$
=\sum_{\vert\alpha\vert\geq 0}(-i)^{\vert\alpha\vert} \sum_{j=1}^n (i(\Psi_{j, \alpha}, (\frac {\partial} {\partial\xi_j} D^{\alpha}f)) + (\Psi_{j,(\alpha_1,\ldots,\alpha_{j-1},\ldots,\alpha_n)},D^{\alpha} f), 
$$
For arbitrary $N \in \mathbb N$ and $j=1,\ldots,n$ define sets
$$
B_N = \{\alpha=(\alpha_1,\ldots,\alpha_n)\in \mathbb Z^n: \alpha_1\leq N,\ldots,\alpha_n\leq N \},
$$
$$
R_{N, j} = \{\alpha_1\leq N,\ldots,\alpha_j = N,\ldots,\alpha_n\leq N, \alpha\in\mathbb Z_+^n\}
$$
and a functional $T_N$ on $G_M(U)$ by the rule
$$
(T_N, f) = \sum_{\alpha \in B_N} (-i)^{\vert\alpha\vert} 
\sum_{j=1}^n (i(\Psi_{j, \alpha}, (\frac {\partial} {\partial\xi_j} D^{\alpha}f)(\xi)) + (\Psi_{j,(\alpha_1,\ldots,\alpha_{j-1},\ldots,\alpha_n)},D^{\alpha} f).
$$
Then $(T, f) = \displaystyle\lim_{N\rightarrow\infty} (T_N, f),\ \ f\in G_M(U)$. 

From the representation
$$
(T_N, f) = 
\sum_{j=1}^n (\sum_{\alpha\in B_N} ((-i)^{\vert\alpha\vert} i (\Psi_{j, \alpha},\frac{\partial} {\partial\xi_j} D^{\alpha} f) + 
$$
$$
+
\sum_{\beta\in B_N}(-i)^{\vert\beta\vert} (\Psi_{j,(\beta_1,\ldots,\beta_{j-1},\ldots,\beta_n)}, D^{\beta} f)). 
$$
we see that for fixed $j \in \{1, \ldots, n \}$ the terms corresponding to a multiindex $\alpha$ with $\alpha_1 \leq N, \ldots,\alpha_j \leq N-1, \ldots,\alpha_n \leq N$, and the summands corresponding to a multiindex
$\beta = (\beta_1, \ldots, \ldots, \beta_n)= \alpha_n$ with $\beta_1 = \alpha_1, \ldots, \beta_j = \alpha_j+1, \ldots, \beta_n= \alpha_n$ annihilate each other. 
Hence
$$
(T_N, f) = \sum_{j=1}^n \sum_{\alpha\in R_{N, j}}(-i)^{\vert\alpha\vert} 
i (\Psi_{j, \alpha},\frac {\partial} {\partial \xi_j} D^{\alpha}f), \ f \in G_M(U).
$$
Now, using (31)  we have for $f\in G_M(U)$
$$
\vert (T_N, f) \vert \leq\sum_{j=1}^n \sum_{\alpha \in R_{N, j}} \vert (\Psi_{j, \alpha},\frac {\partial} {\partial \xi_j} D^{\alpha}f)  \vert 
\leq   
$$
$$
\leq \sum_{j=1}^n \sum_{\alpha \in R_{N, j}} \frac {d_4} {\varepsilon_k^{\vert\alpha\vert} M_{\vert\alpha\vert}} 
\sup_{\xi\in U,\vert\gamma\vert\leq p} (\vert D^{\gamma}(\frac {\partial} {\partial \xi_j}D^{\alpha} f)(\xi) \vert (1 + \Vert \xi \Vert)^p)= 
$$ 
$$
= \sum_{j=1}^n \sum_{\alpha \in R_{N, j}} 
\frac {d_4} {\varepsilon_k^{\vert\alpha\vert} M_{\vert\alpha\vert}} \sup_{\xi\in U, \vert \gamma \vert \leq p}(\vert(D^{(\alpha_1+\gamma_1, \ldots\alpha_j+\gamma_j+1,\ldots,\alpha_n + \gamma_n)} f)(\xi)\vert (1 + \Vert \xi \Vert)^p)).
$$ 
Choose an integer $s >k$ so that $q=\frac{\varepsilon_sH_2} {\varepsilon_k}  < 1 $. 
Then for $f \in G_M(U)$ 
$$
\vert (T_N, f) \vert \leq \sum_{j=1}^n \sum_{\alpha \in R_{N, j}} \frac {d_4} {\varepsilon_k^{\vert\alpha\vert} M_{\vert\alpha\vert}} p_s (f)
\sup_{\xi \in U, \vert\gamma\vert\leq p} \frac 
{\varepsilon_s^{\vert\alpha\vert+\vert\gamma\vert+1} M_{\vert\alpha\vert+\vert\gamma\vert+1}}  {(1 + \Vert\xi\Vert)^{s-p}} \leq 
$$ 
$$
\leq \sum_{j=1}^n \sum_{\alpha \in R_{N, j}} \frac {d_4} {\varepsilon_k^{\vert\alpha\vert} M_{\vert\alpha\vert}} p_s (f) 
\varepsilon_s^{\vert\alpha\vert }\sup_{\vert\gamma\vert\leq p} \varepsilon_s^{\vert\gamma\vert +1 \vert} H_1 H_2^{\vert\alpha\vert+\vert\gamma\vert +1} M_{\vert\alpha\vert} M_{\vert\gamma\vert +1} \leq
$$ 
$$
\leq 
\sum_{j=1}^n 
\sum_{\alpha \in R_{N, j}} d_4 M_{p +1}  
p_s (f) \left(\frac{\varepsilon_s H_2} {\varepsilon_k} \right)^{\vert\alpha\vert} <   
$$ 
$$
< d_4 M_{p +1} p_s (f) n q^N (N+1)^{n-1}.
$$
From this it follows that  for each $f\in G_M(U)$ \ $(T_N, f) \rightarrow 0$ as $N \rightarrow \infty$. 
Thus, $(T, f) =0, \ f \in G_M(U)$. 
So, $T$ is a zero functional. We have proved that the map $L$ is one-to-one.

By the open mapping theorem [17]  $L^{-1}$ is continuous. 
Thus, $L$ is a topological isomorphism. 

The proof of theorem 3 is complete.

\begin{center}
{\bf \S 6. On a dual space for $E(U)$} 
\end{center}

Let $S \in G_M'(U)$. Define a functional $T$ on $E(U)$ by the formula 
$$
(T, F) = (S, f), \ F \in E(U), \ f=F_|U.
$$
Obviously, $T$ is a linear continuous functional on $E(U)$. Let the map ${\cal B}$ acts from $G_M'(U)$ to $E'(U)$ by the rule ${\cal B}(S)= T$.
It is easy to see that the mapping $\cal B$ establishes topological isomorphism of spaces $S\in G_M^{*}(U)$ and $E^{*}(U)$. 

From this and the theorem 3 we obtain the following 

{\bf Theorem 4.} {\it The Laplace transform establishes topological isomorphism between the  spaces $E^*(U)$ and $H_M(T_C)$}.

\pagebreak

I.Kh. Musin, 

Institute of mathematics with Computer centre, Chernyshevskii str., 112, 

Ufa, 450077, Russia

E-mail address: musin@matem.anrb.ru  

P.V. Fedotova,

Bashkirian state university, Ufa, 450000, Rusia

E-mail address: polina81@rambler.ru


\begin{thebibliography}{99}

\bibitem{1} V.S. Vladimirov, {\it Generalized functions in mathematical physics}, Nauka, M.,  1979. 

\bibitem{2} J.W. de. Roever, {\it  Complex Fourier transformation and analytic functionals with unbounded carriers}, Mathematisch Centrum, Amsterdam, 1977. 

\bibitem{3} I. Kh. Musin, P.V. Fedotova,  A theorem of Paley-Wiener type for ultradistributions, {\it Matematicheskie Zametki}, 2009, {\bf 85}:6, 894-914 (in russian), {\it Mathematical Notes}, {\bf 85}:6, 848-867.

\bibitem{4} B.A. Taylor,  Analytically uniform spaces of infinitely differentiable functions, 
{\it Communications on pure and applied mathematics}, {\bf 24}:1 (1971), 39-51.

\bibitem{5} V.S. Vladimirov, Yu. N. Drozhzhinov, B.I. Zav'yalov, {\it Multidimensional tauber theorems for  generalized functions}, Nauks, M., 1986.

\bibitem{6} V.S. Vladimirov,  Functions holomorphic in tubular domains, {\it Izv. Akad. Nauk SSSR Ser. Mat.}, {\bf 27}:1 (1963), 75-100. (Russian)

\bibitem{7} V.S. Vladimirov,  {\it Methods of the theory of functions of many complex variables},  Nauka, M.,  1964.

\bibitem{8} J.W. de Roever, Analytic representation and Fourier transforms of analytic functional in $Z'$ carried by the real space, 
{\it  SIAM J. Math. Anal}, {\bf 9}:6 (1978), 996-1019.

\bibitem{9} J. Sebasti\`ao e Silva,  Su certe classi di spazi localmente convessi importanti per le applicazioni, {\it  Rend. Mat. e Appl.} 
{\bf 14} (1955).

\bibitem{10} V.V. Zharinov,  Compact families of locally convex topological vector spaces, Frechet-Schwartz and dual Frechet-Schwartz spaces, {\it Uspekhi Mat. Nauk }, {\bf 34}:4 (1979), 97-131.

\bibitem{11} L. H\"ormander, {\it An introduction to complex analysis in several variables}, D. Van Nostrand Company,  1966.

\bibitem{12} L. H\"ormander, $L^2$ estimates and existence theorems for the $\overline \partial $ operator, {\it Acta Mathematica}, {\bf 113} (1965), no 1-2, 89-152.

\bibitem{13} V.V. Napalkov,  {\it Convolution equations in multidimensional spaces},  Nauka, M., 1982.

\bibitem{14} R.S. Yulmukhametov, Entire functions of several variables with given behaviour at infinity, 
{\it  Izvestiya Mathematics}, {\bf 60}:4 (1996), 205-224.

\bibitem{15} R.E. Edwards, {\it Functional analysis}, Holt, Pineart and Winston, 1965.

\bibitem{16} S. Mandelbrojt, {\it S\'eries adh\'erentes, r\'egularisation des suites, applications}, Gauthier-Villars, Paris, 1952A.

\bibitem{17} A.P. Robertson and W. Robertson, {\it Topological vector spaces}, Cambridge Univ. Press, Cambridge, 1980.

\end{thebibliography}
\end{document}